\documentclass[10pt,twocolumn,twoside]{IEEEtran}  

\IEEEoverridecommandlockouts                              
\usepackage{amssymb}
\usepackage{amsfonts}
\usepackage{latexsym, amssymb,amsmath, amsbsy,amsopn, amstext}
\usepackage{graphicx}
\usepackage{tikz}
\usepackage{amsmath, amsbsy}
\usepackage{hyperref}
\usepackage{cancel, color}
\usepackage{cite}
\usepackage{float}
\usepackage{ifpdf}
\usepackage{xcolor}
\usepackage{tikz}
\usetikzlibrary{arrows,shapes,snakes,shadows,positioning,automata,patterns}
\usetikzlibrary{trees,decorations.pathmorphing,decorations.markings}

\def\diag{diag}

\def\0{{\bf 0}}

\newtheorem{thm}{Theorem}
\newtheorem{lem}{Lemma}

\newtheorem{rem}{Remark}
\newtheorem{alg}{Algorithm}
\newtheorem{assum}{Assumption}

\title{\LARGE \bf
Distributed generalized Nash equilibria computation  of  monotone  games via preconditioned proximal point algorithms
}

\author{Peng Yi  and Lacra Pavel
\thanks{This work was supported by NSERC Discovery Grant (261764).}
\thanks{P. Yi and L. Pavel are with Department of Electrical and Computer Engineering, University of Toronto, Canada.
        {\tt\small peng.yi@utoronto.ca,pavel@control.toronto.edu}}%
}

\begin{document}

\maketitle
\thispagestyle{empty}
\pagestyle{empty}

\begin{abstract}
In this paper, we investigate distributed  generalized Nash equilibrium (GNE) computation of {\it monotone} games with affine coupling constraints.  Each player can only utilize its local objective function, local feasible set and
a local block of the coupling constraint, and can only communicate with its neighbours.  We assume the game has {\it monotone pseudo-subdifferential} without Lipschitz continuity restrictions. We design novel center-free distributed GNE seeking algorithms for equality and inequality  affine coupling constraints, respectively. A proximal {\bf {\it alternating direction method of multipliers}} (ADMM) is proposed for the equality case, while for the inequality case, a parallel splitting type algorithm is proposed. In both algorithms, the GNE seeking task is decomposed into a sequential NE computation of regularized subgames and distributed update of multipliers and auxiliary variables, based on local data and local communication. { Our two double-layer GNE algorithms need not specify the inner-loop NE seeking algorithm and moreover, only require  that the {\it strongly monotone} subgames are inexactly solved.} {  We prove their convergence by showing that the two algorithms can be seen as specific instances of  {\it preconditioned proximal point algorithms} (PPPA) for finding zeros of monotone operators.}
Applications and numerical simulations are given for illustration.
\end{abstract}

\section{Introductions }\label{sec_introduction}

Generalized Nash equilibrium and its distributed computation is an important research topic in decision making problems over large-scale multi-agent networks. Examples include
power allocation over cognitive radio networks, \cite{pang,wangjian,scutari}, demand response and electric vehicle charging management in smart grids, \cite{hu,zhuminghui,grammatico_2,lygeros2},
rate control over optical networks, \cite{pavel1,pavel2},  and opinion evolution over social networks, \cite{opion,lygeros1}.
Each agent (player)  controls its  decision, and has an objective function to be optimized, which depends on other players' decisions. Moreover, each player's feasible  set can depend on  other players' decisions through coupling constraints, such as { when they share limited network resources.}
Generalized Nash Equilibrium (GNE), firstly proposed in \cite{debreu}, is a reasonable solution,
since at a GNE no player can  decrease/increase its  cost/utility by unilaterally changing its local decision to another feasible one. Interested readers can refer to \cite{faccinei2} for a review on GNE.

Distributed GNE computation methods are quite appealing for noncooperative  games over large-scale networks, in which the local data of each player, including own objective function and own feasible set, are kept by each player. Moreover, when the coupling constraint is a sum of separable local functions, it is also appealing to have each player only knowing its local constraint function, i.e., local contribution to the coupling constraint.
{ Since local data is not required to be transmitted to a central node, the communication burden could be relieved,  and the privacy of each player gets protected.}
Recently,  distributed NE/GNE computation methods have received increasing research attention, see \cite{scutari,zhuminghui,grammatico_2,wangjian,hu,lygeros2} and \cite{baser,yipeng,liangshu,sayed,pavel3,pavel4,shanbhag,shanbhag1,shanbhag4,grammatico_1,lou}.
Different information structures  are considered, depending on whether or not there exists a coordination center.
For example, the methods in \cite{wangjian}\cite{grammatico_2}\cite{lygeros2} all utilize a central node to update and broadcast certain coordination/incentive signals based on all players' decisions. {  Notice that \cite{grammatico_2} considers aggregative games where the agents are coupled through  aggregative variables, hence,  it is efficient to adopt a coordination center if permitted.}
Meanwhile, totally center-free distributed GNE computation algorithms have been proposed in \cite{zhuminghui,sayed,yipeng,liangshu} assuming that each player is able to observe the decisions on which its local objective function or constraint function {\it explicitly} depends on. On the other hand, in the distributed NE computation algorithms of \cite{shanbhag1,pavel3,pavel4}, each player is only required to have local communications with its neighbours, and each player computes an {\it estimation} of other players' decisions or aggregative variables by resorting to consensus dynamics.

Typically, the objective function of each player is convex only with respect to its own decision.
Then an NE/GNE can be computed by solving a (generalized) {\it Variational Inequality} (VI) problem  constructed with  the game's pseudo-gradient/subdifferential (PG/PS) \cite{pang,scutari,wangjian,zhuminghui,faccinei2}.
Various monotonicity and Lipschitz continuity assumptions on PG/PS play a fundamental role in the design and analysis of distributed NE/GNE seeking algorithms.
{ \cite{wangjian} assumes a strongly monotone PG to get the cocoercivity of
the dual operator, and show the convergence of {\it double-layer} dual gradient GNE seeking methods.}
 \cite{sayed} and \cite{yipeng} combine strong monotonicity  and Lipschitz continuity to ensure  the cocoercivity of PG, and propose primal/primal-dual gradient methods for distributed GNE computation.
 \cite{lygeros2,lygeros1} and \cite{grammatico_1} consider {\it aggregative games} with quadratic objective functions, hence also adopt a strong monotone and Lipschitz PG.  \cite{pavel3} and \cite{liangshu} consider games with strictly monotone and Lipschitz PG. \cite{pavel3} proposes a ``gradient"+``consensus" algorithm for distributed NE seeking, while \cite{liangshu} utilizes a continuous-time gradient flow algorithm to seek a GNE of aggregative games.
{
For NE seeking with only {\it monotone PGs},
\cite{scutari} proposes a {\it double-layer} proximal best-response algorithm that involves solving regularized subgames at each iteration, while \cite{shanbhag} proposes a  single time-scale/layer regularized (sub)gradient algorithm with {\it diminishing step-sizes}. For  GNE seeking of {\it monotone games},
\cite{pavel2} proposes a double-layer dual extragradient method and \cite{shanbhag4} adopts the  single-layer Tikhonov regularization algorithm with {\it diminishing step-sizes}, both assuming Lipschitz continuity and using a central coordinator.
\cite{zhuminghui} proposes a primal-dual gradient algorithm, and \cite{payoff} proposes a payoff-based algorithm for  GNE seeking with  pseudo-monotone  PGs, both with diminishing step-sizes.
}

Motivated by the above, we investigate center-free distributed algorithms for computing GNE of {\it monotone games with affine coupling constraints}.  The players' decisions are coupled together with a globally shared affine constraint, while each player only knows a local block of the constraint. We consider both equality and inequality constraints which cover many task/resource allocation games, \cite{wangjian,lygeros2,liangshu,yipeng}.
{\it Compared with previous works, the key difference is that we only assume a monotone pseudo-subdifferential without Lipschitz continuity restrictions.} We  propose {\it center-free GNE algorithms with fixed step-sizes} where each player only utilizes its local data and  has a peer-to-peer communication with its neighbours. 
To the best of our knowledge, this distributed GNE computation has not been discussed in literature under this general form.

We adopt the variational GNE  as a refined solution and  use primal-dual analysis to reformulate  GNE seeking  as the problem of finding  zeros of  monotone operators for equality and inequality cases, respectively. The monotone operators are composed of a skew-symmetric linear operator (with both the constraint matrices and a matrix related to communication graph) and an operator involving  PS.
In general, the proximal point algorithm can be applied for solving  monotone inclusion problems without Lipschitz restriction.
However, it is not directly applicable to our GNE problem
because it requires to compute the inverse of a graph-related skew-symmetric matrix, which is prohibitive in distributed algorithms. To overcome these challenges, we propose novel distributed GNE seeking algorithms based on Preconditioned Proximal Point Algorithm (PPPA), for equality and inequality cases, respectively.
For the equality  case, we call it proximal {\it alternating direction method of multipliers} (ADMM), partially motivated by  \cite{hebingsheng}. For the inequality case, we call it  proximal parallel splitting algorithm, partially motivated by \cite{hebingsheng2}.  Both algorithms use appropriately chosen operators and preconditioning matrices, which ensure that the resolvent evaluation of monotone operators is realizable  by local computation and communication.
The proposed algorithms decompose the GNE computation  into sequential  NE computation for regularized subgames and distributed update of local multipliers and auxiliary variables. { Hence, our algorithms are double-layer algorithms, similar to \cite{scutari,wangjian,pavel2}, but the inner-loop NE seeking algorithms need not be specified while the subgame only needs to be solved inexactly.} By using proximal terms, the subgame is regularized to have strongly monotone PS, hence it can be efficiently solved by existing NE seeking distributed algorithms, such as the  best-response algorithm in \cite{scutari}. { The inexactness in solving the subgames is also considered and relaxation steps are applied to all variables, which potentially could improve convergence speed.} 
 Moreover, proximal ADMM enjoys the feature of utilizing the most recent available information whenever possible. In both cases, the algorithms'  convergence is proved for {\it fixed step-sizes} by relating them to PPPA, and showing that they can be seen as specific instances of PPPA, while PPPA's convergence can be shown based on averaged operator theory.

{
To summarize, the main contributions of this work are as follows.
(i): The game model only assumes a {\it monotone PS without Lipschitz continuity}, hence it is a generalization of previous ones. Both equality and inequality affine coupling constraints are considered.
(ii): Novel center-free GNE seeking algorithms with peer-to-peer communication are introduced.
Since only monotonicity is imposed, the double-layer algorithms could be implemented after the NE algorithm is chosen tailored to the specific practical problem.
Moreover, thanks to  the proximal terms, the subgames are regularized to have strongly monotone PS/PGs, hence, could be efficiently solved.
(iii): The algorithms are related to  PPPA for monotone inclusion, revealing the algorithms' intrinsic structure.
Their convergence is proved for fixed step-sizes.
}

The paper is organized as follows.
Section \ref{sec_notations_and preliminaries} gives the preliminary background.
Section \ref{sec_game_formulation} formulates the noncooperative game and  basic assumptions.
Section \ref{sec_algorithm_and_limiting_poinit} gives distributed GNE computation algorithms for both equality and inequality constraint cases,  and analyzes their limiting points.
Section \ref{sec_convergence_analysis} presents the algorithms' convergence analysis.
Section \ref{sec_simulations} gives application examples and simulation studies.
Section \ref{sec_concluding} draws the  concluding remarks.

\section{notations and preliminaries}\label{sec_notations_and preliminaries}
In this section, we review the notations and preliminary notions in  monotone and averaged operators from \cite{combettes1}.

{\it Notations}:  In the following,
$\mathbf{R}^m$ ($\mathbf{R}^m_{+}$) denotes the $m-$dimesional (nonnegative) Euclidean space.
For a column vector $x \in \mathbf{R}^m$ (matrix $A\in \mathbf{R}^{m\times n}$),
$x^T$ ($A^T$)  denotes its transpose.
$x^Ty=\langle x, y\rangle$ denotes the inner product of $x,y$, and $||x||= \sqrt{x^Tx}$ denotes the induced norm. 
$||x||^2_G$ denotes $\langle x, Gx\rangle$ for a symmetric matrix $G$.
Denote $\mathbf{1}_m=(1,...,1)^T \in \mathbf{R}^m$ and
$\mathbf{0}_m=(0,...,0)^T \in \mathbf{R}^m$.
$diag \{A_1, . . . ,A_N\}$ represents
the block diagonal matrix with  $A_1, . . . ,A_N$ on its
main diagonal. Denote $col(x_1,....,x_N) $ as the stacked column vector of $x_1$ to $x_N$.
$I_n$ denotes the identity matrix in $\mathbf{R}^{n\times n}$.
For a matrix $A=[a_{ij}]$, $a_{ij}$ or $[A]_{ij}$
stands for the matrix entry in the $i$th row and $j$th column of $A$.
Denote $int(\Omega)$ as the interior of $\Omega$ and
$ri(\Omega)$ as the relative interior of $\Omega$.
Denote $\times_{i=1,...,N}\Omega_i$ or $\prod_{i=1}^N \Omega_i$ as the Cartesian product of  $\Omega_i,i=1,...,N$.

Let $\mathfrak{A}:\mathbf{R}^m \rightarrow 2^{\mathbf{R}^m}$ be a set-valued operator. ${\rm Id}$ denotes the identity operator, i.e, ${\rm Id}(x)=x$.
The domain of $\mathfrak{A}$ is $dom\mathfrak{A}= \{x\in \mathbf{R}^m| \mathfrak{A}x \neq \emptyset\}$ where $\emptyset$ stands for the empty set, and the range of $\mathfrak{A}$ is $ran\mathfrak{A}=\{y \in \mathbf{R}^m| \exists x, y\in \mathfrak{A}x\}$. The graph of $\mathfrak{A}$ is $gra\mathfrak{A}=\{(x,u) \in \mathbf{R}^m\times \mathbf{R}^m| u\in \mathfrak{A}x\}$. The inverse of $\mathfrak{A}$ is defined via $gra\mathfrak{A}^{-1}=\{(u, x)| (x, u)\in gra \mathfrak{A}\}$.
The zero set of  $\mathfrak{A}$ is $zer\mathfrak{A}=\{x\in \mathbf{R}^m | \mathbf{0} \in \mathfrak{A}x\}$.
The sum of  $\mathfrak{A}$ and $\mathfrak{B}$ is defined as $gra (\mathfrak{A}+\mathfrak{B})=\{(x,y+z)| (x,y)\in gra \mathfrak{A}, (x,z)\in gra \mathfrak{B}\}$. Define the {\it resolvent} of $\mathfrak{A}$ as $R_{\mathfrak{A}}=({\rm Id}+\mathfrak{A})^{-1}$.

Operator $\mathfrak{A}$ is monotone if
$\forall (x,u), \forall(y,v)\in gra\mathfrak{A}$, we have
$\langle x-y, u-v\rangle \geq 0.$
 $\mathfrak{A}$ is maximally monotone if $gra\mathfrak{A}$ is not {\it strictly} contained in the graph of any other monotone operator.
A skew-symmetric matrix $A=-A^T$ defines a maximally monotone operator $Ax$ (\cite{combettes1}, p. $298$).
Suppose  $\mathfrak{A}$ and $\mathfrak{B}$ are maximally monotone operators and
$0\in int(dom \mathfrak{A}-dom \mathfrak{B})$, then $\mathfrak{A}+\mathfrak{B}$ is also maximally monotone.
For a proper {\it lower semi-continuous convex} (l.s.c.) function $f$, its subdifferential operator $\partial f: domf\rightarrow 2^{\mathbf{R}^m}$ is
$
\partial f:\; x \mapsto \{g | f(y)\geq f(x)+ \langle g, y-x \rangle, \forall y\in domf\}.
$
$\partial f$ is  maximally monotone and $Prox_{f}= R_{\partial f}:\mathbf{R}^m\rightarrow dom f$ is called the proximal operator of $f$,
i.e., $ Prox_{f} :  x \mapsto  \arg\min_{u\in dom f } f(u)+\frac{1}{2} ||u-x ||_2^2.$

Define the indicator function of  $\Omega$ as $\iota_{\Omega}(x)= 0$ if $ x\in \Omega$  and $\iota_{\Omega}(x)= \infty$ if $x\notin \Omega.$
For a closed convex set $\Omega$, $\iota_{\Omega}$ is a proper l.s.c. function.
$\partial \iota_{\Omega}$ is also the normal cone operator of  $\Omega$, i.e., $N_{\Omega}(x)$,  where
   $N_{\Omega}(x)=\{v| \langle v, y-x\rangle\leq 0, \forall y\in \Omega\}  $ and $dom N_{\Omega}=\Omega$.
Given a symmetric positive definite matrix $G$, define $P^G_{\Omega}(x)=\arg\min_{y} (\iota_{\Omega}(y)+ \frac{1}{2}||x-y ||_G^2)$.

For a single-valued operator $T:\Omega \subset \mathbf{R}^m\rightarrow \mathbf{R}^m$,  $x\in \Omega$ is a fixed point of $T$ if $Tx=x$. 
$T$  is nonexpansive if
it is $1-$Lipschitzian, i.e., $||T(x)-T(y) || \leq ||x-y||, \forall x,y \in \Omega$. $T$ is contractive if
$\exists \gamma \in (0,1)$ s.t. $ ||T(x)-T(y) || \leq \gamma ||x-y||, \forall x,y \in \Omega.$
Let $\alpha \in (0,1)$, then $T$ is $\alpha-$averaged, denoted as  $T\in \mathcal{A}(\alpha) $,  if $\exists$ a nonexpansive operator $T^{'}$ such that  $T=(1-\alpha){\rm Id}+\alpha T^{'}$.
If $T\in \mathcal{A}(\frac{1}{2})$,  $T$ is called firmly nonexpansive.

\section{Game formulation}\label{sec_game_formulation}
Consider a set of players (agents) $\mathcal{N}=\{1,\cdots, N \}$ that are involved in the following noncooperative game with shared coupling constraints.
Player $i\in \mathcal{N}$ controls its own decision (strategy or action) $x_i\in \Omega_i \subset \mathbf{R}^{n_i}$, where $\Omega_i $ is its private feasible set.  Let  $\mathbf{x} =col(x_1,\cdots,x_N) \in \mathbf{R}^n$ denote the decision profile,  i.e., the stacked vector of all agents' decisions, with $\sum_{i=1}^N n_i=n$.
Let  $\mathbf{x}_{-i}=col(x_1,\cdots,x_{i-1},x_{i+1},\cdots,x_{N})$ denote the decision profile of all agents except player $i$. Player $i$ aims to optimize its own objective function within its feasible set,  $f_i(x_i,\mathbf{x}_{-i}): \bar{\Omega} \rightarrow \mathbf{R}$ where $\bar{\Omega}=\prod_{i=1}^N \Omega_i  \subset \mathbf{R}^{n}$. Note that $f_i(x_i,\mathbf{x}_{-i})$ is coupled with  other players' decisions $\mathbf{x}_{-i}$.
Moreover, all the players' decisions are coupled together through a {\it globally shared set } $X \subset \mathbf{R}^n$.
Hence, player $i$ has a set-valued map $X_i(\mathbf{x}_{-i}): \mathbf{R}^{n-n_i}\rightarrow 2^{\mathbf{R}^{n_i}}$ that specifies its feasible set  defined as
$$X_i(\mathbf{x}_{-i}): =\{x_i\in  \Omega_i| (x_i,\mathbf{x}_{-i})\in X\}.$$ 
Given $\mathbf{x}_{-i}$, player $i$'s best-response strategy is
\begin{equation}\label{GM}
\min_{x_i} \;  f_i(x_i,\mathbf{x}_{-i}),  \; s.t., \; \; x_i \in X_i(\mathbf{x}_{-i}).
\end{equation}
A generalized Nash equilibrium (GNE) $\mathbf{x}^*=col(x_1^*,\cdots,x_N^*)$ is defined at the intersection of all  players' best-response sets,
\begin{equation}\label{GNE}
x_i^*\in \arg\min_{x_i} f_i(x_i,\mathbf{x}^*_{-i}), \; s.t., \;  x_i \in X_i(\mathbf{x}^*_{-i}), \; \forall i\in \mathcal{N}.
\end{equation}

We consider the set $X$ defined via two types of shared affine coupling constraints, equality and inequality constraints.
For the equality constraint case, $X=X^e$ where we denote
\begin{equation}\label{coupling_set_equ}
 X^e  := \prod_{i=1}^N \Omega_i \bigcap \{ \mathbf{x} \in \mathbf{R}^n |  \sum_{i=1}^N A_ix_i= \sum_{i=1}^N b_i\}. 
\end{equation}
For the inequality constraint case, $X=X^i$ where
\begin{equation}\label{coupling_set_ine}
{X}^i := \prod_{i=1}^N \Omega_i \bigcap \{ \mathbf{x} \in \mathbf{R}^n |  \sum_{i=1}^N A_ix_i \leq  \sum_{i=1}^N b_i\}.
\end{equation}
In both \eqref{coupling_set_equ} and \eqref{coupling_set_ine}, $A_i\in \mathbf{R}^{m\times n_i}$ and $b_i\in \mathbf{R}^m$ as well as $\Omega_i$ are private data of player $i$. 
 Thereby, the shared  set $X$ couples all players' feasible sets, but is not known by any agent.
We consider the following  assumption on the game in \eqref{GM}.

{ {\it \begin{assum}\label{assum1}
For player $i$, $f_i(x_i,\mathbf{x}_{-i})$ is a proper l.s.c. function with respect to $x_i$ given any fixed $\mathbf{x}_{-i}$, and its  subdifferential with respect to $x_i$ is $\partial_i f_i(x_i,\mathbf{x}_{-i}) $. The {\it pseudo-subdifferential} of the game in \eqref{GM} defined as  $\partial F(\mathbf{x}): \mathbf{x}\rightarrow \prod_{i=1}^N\partial_i f_i(x_i,\mathbf{x}_{-i})$  is maximally monotone.
$\Omega_i$ is  a closed  convex set with nonempty interior.  $X^e$ in \eqref{coupling_set_equ} has nonempty relative interiors, and $X^i$ in \eqref{coupling_set_ine} has nonempty interiors.
$X_i(\mathbf{x}_{-i})$ has nonempty relative interiors for $\mathbf{x}_{-i} \in \prod_{j=1, j \neq i}^N \Omega_j$ when $X=X^e$, and
$X_i(\mathbf{x}_{-i})$ has nonempty  interiors for $\mathbf{x}_{-i} \in \prod_{j=1, j \neq i}^N \Omega_j$ when $X=X^i$.
\end{assum}}}

\begin{rem}\label{rem_splitting}
In many practical cases, $f_i(x_i,\mathbf{x}_{-i})$  has a splitting structure such as $f_i(x_i,\mathbf{x}_{-i})= g_i(x_i,\mathbf{x}_{-i})+l_i(x_i)$, \cite{shanbhag4}, where $g_i(x_i,\mathbf{x}_{-i})$ is differentiable and convex with respect to $x_i$, and 
$l_i(x_i)$ is a local l.s.c. regularization/cost term.
Denote
$\nabla_p G(\mathbf{x})=col(\nabla_{1} g_1(x_1,\mathbf{x}_{-1}),\cdots, \nabla_{N} g_N(x_N,\mathbf{x}_{-N}) )$ where $\nabla_{i} g_i(x_i,\mathbf{x}_{-i})$ is the gradient of $g_i$ with respect to $x_i$ and
 $\partial L(\mathbf{x}): \mathbf{x}\rightarrow \prod_{i=1}^N \partial l_i(x_i)$. Then $\partial L(\mathbf{x})$ is  maximally monotone,  since it is the subdifferential of $\sum_{i=1}^N l_i(x_i)$.
In this case, $\partial F(\mathbf{x})=\partial L(\mathbf{x})+ \nabla_p G(\mathbf{x})$ is maximally monotone when $\nabla_p G(\mathbf{x})$ is monotone.
\end{rem}

Define the {\it generalized variational inequality} (GVI) problem 
\begin{equation}\label{vi_original}
Find \; \mathbf{x}^*, \; s.t. \; \langle l^*,\mathbf{x}-\mathbf{x}^*\rangle \geq 0, l^*\in \partial F(\mathbf{x}^*), \forall \mathbf{x}\in  X.
\end{equation}
According to Proposition 12.4 in \cite{pang}, any solution of \eqref{vi_original} is a GNE of game in \eqref{GM}, called {\it variational GNE}.

Let us first analyze the equality constraint case,  $X=X^e$.
Under Assumption \ref{assum1}, $\mathbf{x}^*$ is a  GNE of the game in \eqref{GM} if and only if
$\forall i\in \mathcal{N}$ there exists $\lambda_i^* \in \mathbf{R}^m$ such that,
 \begin{equation}
 \begin{array}{l}\label{kkt1}
 \mathbf{0}\in \partial_i f_i(x_i^*,\mathbf{x}^*_{-i})+A_i^T \lambda_i^*+N_{\Omega_i}(x_i^*),\; \forall i\in \mathcal{N},\\
\sum_{i=1}^N A_i x_i^* =\sum_{i=1}^N b_i.
\end{array}
 \end{equation}
 Meanwhile, based on the Lagrangian duality for GVI (Equation (12.4) of \cite{pang}), $\mathbf{x}^*$ is a solution of  GVI in \eqref{vi_original} with $X=X^e$ if and only if there exists a multiplier $\lambda^* \in \mathbf{R}^m $ such that 
 \begin{equation}
 \begin{array}{l}\label{kkt2}
  \mathbf{0} \in \partial_i f_i(x_i^*,\mathbf{x}^*_{-i}) + A_i^T \lambda^*+N_{\Omega_i}(x_i^*), \quad \forall i\in \mathcal{N},\\
\sum_{i=1}^N A_i x_i^* =\sum_{i=1}^N b_i.
 \end{array}
 \end{equation}
 By comparing the KKT conditions in \eqref{kkt1} and \eqref{kkt2}, we have that any solution to GVI in \eqref{vi_original} with $X=X^e$ is a GNE of the game in \eqref{GM} with all players having the same local  multiplier.

Similarly, for the inequality case  $X=X^i$,
$\mathbf{x}^*$ is a solution of  GVI \eqref{vi_original} with $X=X^i$ if and only if there exists a multiplier $\lambda^* \in \mathbf{R}_{+}^m $ such that 
 \begin{equation}
 \begin{array}{l}\label{kkt2_ine}
  \mathbf{0} \in \partial_i f_i(x_i^*,\mathbf{x}^*_{-i}) + A_i^T \lambda^*+N_{\Omega_i}(x_i^*), \quad \forall i\in \mathcal{N},\\
  \mathbf{0}\in -\sum_{i=1}^N (A_ix_i^*-b_i)+ N_{\mathbf{R}^m_{+}}(\lambda^*).
 \end{array}
 \end{equation}

Not every GNE of the considered game in \eqref{GM} is a solution to the GVI in \eqref{vi_original}.
Since the variational GNE has an economic interpretation of no price discrimination and enjoys a stability and sensitivity property (refer to \cite{pang}),  {\it we aim to propose novel distributed algorithms for computing  a variational GNE of the monotone game} for $X=X^e$ and $X=X^i$, respectively.
{\it \begin{assum}\label{assum2}
The solution set of GVI in \eqref{vi_original} is nonempty for both $X=X^e$ and $X=X^i$, or equivalently, the considered game in \eqref{GM} has at least a variational GNE.
\end{assum}}

\begin{rem}
Some sufficient conditions for the existence of solutions to monotone GVI can be found in \cite{pang} and \cite{shanbhag4}. For example,  compactness of $\Omega_i,\forall i\in \mathcal{N}$ ensures Assumption \ref{assum2}.
\end{rem}

\section{Distributed GNE computation algorithms}\label{sec_algorithm_and_limiting_poinit}

In this section, we propose distributed algorithms that players can use to find a solution of GVI  \eqref{vi_original} for $X=X^e$ and $X=X^i$, respectively.
We focus on distributed variational GNE computation because of two reasons. Firstly,
player $i$ can only manipulate its local  $f_i(x_i,\mathbf{x}_{-i})$, $A_i$, $b_i$ and $\Omega_i$ for local computation, since these contain its private information.
{ Secondly, we assume there is no central node that has bidirectional communications with all players, either because this could be inefficient from a communication point of view, or because it might be not possible to have such a central node. }
Thus, each player only uses  its local data for local computation, and has peer-to-peer communication with its neighbours for local coordination.

We first introduce the communication graph and algorithm notations in  \ref{sub_sec_auxiliary_variables}.
We give the proximal ADMM for equality constraint case in \ref{seb_sec_proximal ADMM},
and the distributed algorithm for inequality constraint case in \ref{sec_parallel_splitting}.

\subsection{Communication graph and algorithm variables}\label{sub_sec_auxiliary_variables}
To facilitate the distributed coordination, players are able to communicate with  their neighbours through a connected and undirected graph $\mathcal{G}=(\mathcal{N},\mathcal{E})$.
The edge set is $\mathcal{E} \subset \mathcal{N}\times \mathcal{N} $, 
$(i,j) \in \mathcal{E}$ if agent $i$ and agent $j$ can share information with each other, and
agents $j$, $i$ are  called neighbours.
A path of graph $\mathcal{G}$ is a sequence of distinct agents in
$\mathcal{N}$ such that any consecutive agents in the sequence
are neighbours. Agent $j$ is said to be
connected to agent $i$ if there is a path from $j$ to $i$.
$\mathcal{G}$ is connected if any two agents are
connected.

Obviously, $|\mathcal{N}|=N$, and we denote $|\mathcal{E}|=M$.
The edges are labeled with $e_l$, $l=1,\cdots,M$.
Without  loss  of  generality, $e_l=(i,j)$ is {\it arbitrarily} ordered and denoted by $i\rightarrow j$.
Define $\mathcal{E}^{in}_i$ and $\mathcal{E}^{out}_i$ for  agent $i$ as follows:
$e_l \in  \mathcal{E}^{in}_i$  if agent $i$ is the targeted point of $e_l$; $e_l\in \mathcal{E}^{out}_i $ if agent $i$ is the starting point of $e_l$. Then denote $\mathcal{E}_i=\mathcal{E}^{in}_i \bigcup \mathcal{E}^{out}_i$ as the set of edges adjoint to agent $i$.
Define the {\it incidence matrix} of $\mathcal{G}$ as ${V}\in \mathbf{R}^{N\times M}$ with
${V}_{il}=1$ if  $e_l\in \mathcal{E}_i^{in}$, and ${V}_{il}=-1$ if $e_l\in \mathcal{E}_i^{out}$, otherwise ${V}_{il}=0$. We have $\mathbf{1}^T_N V=\mathbf{0}^T_M$,
 and $V^Tx=\mathbf{0}_M$ if and only if $x\in \{\alpha \mathbf{1}_N | \alpha \in\mathbf{ R}\}$ when $\mathcal{G}$ is connected.
 Denote $\mathcal{N}_l=\{i,j\}$ as the pair of agents  connected by edge $e_l=(i,j)$.

We introduce the variables.   Firstly, each player has a local decision $x_i\in \Omega_i$ and  a {\it local  multiplier} $\lambda_i\in \mathbf{R}^m$. According to KKT \eqref{kkt2} and \eqref{kkt2_ine}, in steady-state all players should have the same local  multiplier, i.e., $\lambda_i=\lambda^*,\forall i\in \mathcal{N}$.
To facilitate the coordination for the consensus of local multipliers and to ensure the  coupling constraint, we consider an auxiliary variable $z_l\in \mathbf{R}^m$ associated with edge $e_l$ of graph $\mathcal{G}$.
Notice that $\mathcal{G}$ is undirected and  the edges are {\it arbitrarily} ordered, therefore,  we can have any agent from $\mathcal{N}_l$ to  maintain $z_l$.
For clarity, we let the starting agent of an edge to maintain the corresponding edge variable.
That is agent $i$ will take the responsibility for maintaining $z_l$ if $e_l\in  \mathcal{E}^{out}_i$.

Before presenting the algorithms we first make some observations.  The  algorithms are based on decomposing the GNE computation  into sequential NE computation for regularized subgames and distributed update of local multipliers and auxiliary variables. The regularized subgames are made to have strongly monotone PS with the help of proximal terms, hence can be efficiently solved by existing distributed algorithms, such as the  best-response algorithm in \cite{scutari}. The update of the local multipliers has to be done so that in steady-state they are the same, and satisfy the optimality conditions \eqref{kkt2},  \eqref{kkt2_ine}  involving the constraints, while using only local information. Towards this we use the auxiliary variables $z_l$, which have a double role:  to help in estimating the contribution of the other players' in the constraints and to enforce consensus.

Let $x_{i,k}, \lambda_{i,k}$ and $z_{l,k}$ denote  $x_i,\lambda_i$ and $z_l$ at iteration $k$.

\subsection{Proximal ADMM for $X=X^e$}\label{seb_sec_proximal ADMM}

The distributed  algorithm for computing a variational GNE of game in \eqref{GM} when $X=X^e$ is given as follows.

\begin{alg}\label{alg_d}
\quad \\
\noindent\rule{0.49\textwidth}{0.7mm}
{\bf Step 1--update of $x_{i,k}$}:
\begin{itemize}
\item Player $i$ receives $z_{l,k},l\in \mathcal{E}^{in}_i$ through $\mathcal{G}$.

 \item Construct a subgame where player $i$ has a decision $ x_i\in \Omega_i $ and an objective function $\tilde{f}_i(x_i,\mathbf{x}_{-i})$,
\begin{align}\label{subgame}\vspace{-0.2cm}
&\tilde{f}_i(x_i,\mathbf{x}_{-i})=f_i(x_i,\mathbf{x}_{-i})+ \frac{1}{2} || x_i-x_{i,k}||_{R_i}^2\nonumber \\
&\; + [\lambda_{i,k}+H_i(A_ix_{i,k}+\sum_{l\in \mathcal{E}_i} V_{il}z_{l,k}-b_i)]^T A_i x_i.
\end{align}
and denote its  NE by $\hat{\mathbf{x}}_{k}=col(\hat{x}_{1,k},\cdots, \hat{x}_{N,k})$.
\item   Players compute $\tilde{\mathbf{x}}_{k}=col(\tilde{x}_{1,k},\cdots,\tilde{x}_{N,k})$ as an inexact solution to subgame \eqref{subgame} such that  $||\tilde{\mathbf{x}}_{k}-\hat{\mathbf{x}}_{k}|| \leq \mu_{k}$, where $\mu_k$ is described below.

\item  Player $i$ updates its local decision $x_{i,k}$ with
\begin{equation}\label{alg_d_x_dynamics}
x_{i,k+1}=x_{i,k}+\rho (\tilde{x}_{i,k}-x_{i,k}).
\end{equation}
\end{itemize}
\vspace{-0.2cm}

{\bf Step 2--update of $\lambda_{i,k}$}:
\begin{equation}\label{alg_d_lambda_dynamics}
{\lambda}_{i,k+1}=\lambda_{i,k}+\rho H_i(A_i\tilde{x}_{i,k}+\sum_{l\in \mathcal{E}_i} V_{il}z_{l,k}-b_i).
\end{equation}
\vspace{-0.4cm}

{\bf Step 3--update of $z_{l,k}$}:

Let $s_{i,k}= \frac{1}{\rho}\lambda_{i,k+1}+\frac{\rho-1}{\rho}\lambda_{i,k}
+H_i(A_i\tilde{x}_{i,k}+\sum_{l\in \mathcal{E}_i}V_{il}z_{l,k}-b_i).$
For $e_l\in \mathcal{E}_i^{out}$, player $i$ receives $s_{j,k}$, $j\in \mathcal{N}_l\setminus \{i\}$,  and updates $z_{l,k}$ with
\begin{equation}
\begin{array}{l}\label{alg_d_z_dynamics}
z_{l,k+1}=z_{l,k}-\rho W_l (s_{j,k}-s_{i,k}).
\end{array}
\end{equation}
\noindent\rule{0.49\textwidth}{0.7mm}
 $\{\mu_{k}\}$ is a nonnegative sequence s.t. $\sum_{k=1}^{\infty}\mu_{k} <\infty$, $\rho\in [1,2)$ is a fixed relaxation/extrapolation step-size, and $R_i \in \mathbf{R}^{n_i\times n_i}$, $H_i \in \mathbf{R}^{m\times m}$ and $W_l\in \mathbf{R}^{m\times m}$, $l\in \mathcal{E}_i^{out}$ are local parameters (step-sizes) that are symmetric positive definite matrices.
\end{alg}

\vskip 2mm
We give next some  intuition behind Algorithm \ref{alg_d}'s design. Since $A_i$, $b_i$ are private data, the coupling constraint $X=X^e$ is not completely known by any player. Note that in steady-state we should have $A_ix_i^*-b_i=\sum_{j=1,j\neq i}^N (A_jx^*_j-b_j)$ due to \eqref{kkt2}, where the right-hand side is unknown information for player $i$.
The penalized cost  $\tilde{f}_i$ in \eqref{subgame} is composed of  a proximal term $f_i(x_i,x_{-i}) +\frac{1}{2}||x_i-x_{i,k} ||_{R_i}^2$ (to regularize the subgames), a Lagrangian term $ \langle \lambda_i, A_ix_i\rangle $  and a penalty term. The penalty term is based on linearizing the  quadratic penalty $\frac{1}{2}||A_ix_i-b_i+\sum_{j=1,j\neq i}^N (A_j x_{j,k}-b_j)||_{H_i}^2$ at $x_{i,k}$, which should be zero in steady-state cf.  \eqref{kkt2}. This gives $\langle H_i(A_ix_{i,k}+\sum_{j=1,j\neq i}^N (A_j x_{j,k}-b_j)-b_i),A_ix_i\rangle$  after dropping all constants. To overcome the need for information about the other players $j\neq i$,  in \eqref{subgame} this term is estimated  as  $(H_i[A_ix_{i,k}+\sum_{l\in \mathcal{E}_i} V_{il}z_{l,k}-b_i])^TA_ix_i$,  via the auxiliary variables $z_l$. A similar term is used in the local multiplier $\lambda_i$'s update, \eqref{alg_d_lambda_dynamics}. Player $i$ uses $A_i\tilde{x}_{i,k}+\sum_{l\in \mathcal{E}_i} V_{il}z_{l,k}-b_i$ as  an estimation of $\sum_{j=1}^N(A_jx_{j,k}-b_j)$ to update its $\lambda_i$. The update for $z_l$,  \eqref{alg_d_z_dynamics},  has an integrator dynamics form driven by the difference between  $\lambda_i$ and $ \lambda_j$, since $e_l=(i,j)$, and ensures the consensus of local multipliers.
Meanwhile, \eqref{alg_d_lambda_dynamics} also utilizes $z_{l,k}$, that is the integrator  for differences between multipliers, as the feedback signal to reach consensus of local multipliers.

We show in Theorem \ref{thm_limiting_point} that at the limit point of the algorithm  $A_ix_i^*-b_i = \sum_{l\in \mathcal{E}_i} V_{il}z^*_l,\forall i\in \mathcal{N}$, while  $A_ix_i^*-b_i=\sum_{j=1,j\neq i}^N (A_jx^*_j-b_j)$ due to \eqref{kkt2}.  Hence, $\sum_{l\in \mathcal{E}_i} V_{il}z^*_l$, generated as an output of \eqref{alg_d_z_dynamics}  is an  estimation of
$\sum_{j=1,j\neq i}^N (A_jx^*_j-b_j)$, and $\sum_{l\in \mathcal{E}_i} V_{il}z_{l,k}$ as used by player $i$ is a dynamical estimator for $\sum_{j=1,j\neq i}^N (A_jx_{j,k}-b_{j})$.
{  Motivated by \cite{low}, the auxiliary variable $z_l$ has an interpretation of network flow.
In fact, if we regard $A_ix_i$ as in-flow  at node $i$ and $b_i$ as out-flow at node $i$, and  $\sum_{i=1}^N A_i x_i=\sum_{i=1}^N b_i$ is a conservative network flow balancing constraint. Thereby,  $z_l$ can be regarded as  flow on each edge to ensure the balancing constraint.} All in all,  variables $z_l$ estimate  the other players' contribution to coupling constraints, and ensure local multipliers reach consensus.

Algorithm \ref{alg_d} updates each coordinate with the most recent information in a Gauss-Seidel manner and uses proximal terms, hence is called {\it proximal ADMM}.  It uses relaxation steps, $\rho \in [1,2)$, to perform extrapolations of all variables, {  which in practice could accelerate convergence  (refer to Figure 2 of \cite{eckstein} and numerical studies in \cite{hebingsheng,hebingsheng2}).} It  is a center-free distributed algorithm with peer-to-peer communications. In {\bf Step 1}, player $i$ communicates with its neighbours to get $z_{l,k},l\in \mathcal{E}_i^{in}$. The NE of subgames can be computed in a distributed manner with existing algorithms such as  best-response algorithms in \cite{scutari} and gradient algorithms in \cite{shanbhag}, \cite{pavel3} and \cite{pavel4}, which only involve local computations and communications.
In {\bf Step 2}, player $i$ uses its  $\tilde{x}_{i,k}$ and locally available $z_{l,k},l\in \mathcal{E}_i$  to update its local multiplier $\lambda_i$.
In {\bf Step 3}, player $i$ computes $s_{i,k}$ with its local information, and receives $s_{j,k},j\in \mathcal{N}_{l}\setminus \{i\}$ to update $z_{l,k}$.

Next, we put Algorithm \ref{alg_d} in a compact form and  show that  its limiting point $\mathbf{x}^*$ is a variational GNE of game in \eqref{GM} when $X=X^e$.
We use the following compact notations.
Denote $\bar{\lambda}=col(\lambda_1,\cdots,\lambda_N)$ and $\mathbf{Z}=col(z_1,\cdots,z_M)$.
Denote
$R=\diag\{R_{1},\cdots,R_{N}\}$,
$W=\diag\{W_1,\cdots,W_M\}$,
$H=\diag\{H_{1},\cdots, H_{N}\}$,
$\bar{V}=V\otimes I_m$,
$\Lambda=\diag\{A_1,\cdots,A_N\}$,
 and $\bar{b}=col(b_1,\cdots,b_N)$.

\begin{thm}\label{thm_limiting_point}
Suppose that Assumption \ref{assum1} and \ref{assum2} hold for game  \eqref{GM} when $X=X^e$. Then any limiting  point $col(\mathbf{x}^*,\mathbf{Z}^*,\bar{\lambda}^*)$ of Algorithm $\ref{alg_d}$ belongs to the zeros of operator $\mathfrak{M}^e $ defined by\vspace{-0.1cm}
\begin{equation} \label{MI}
\mathfrak{M}^e: \left(
  \begin{array}{c}
    \mathbf{x} \\
    \mathbf{Z} \\
    \bar{\lambda} \\
  \end{array}
\right) \mapsto
\left(
  \begin{array}{c}
   \Lambda^T\bar{\lambda}+  (N_{\bar{\Omega}}+\partial F)\mathbf{x} \\
    \bar{V}^T\bar{\lambda }  \\
    -\Lambda\mathbf{x} -\bar{V}\mathbf{Z}+\bar{b}    \\
  \end{array}
\right)
\end{equation}
Meanwhile, any zero $col(\mathbf{x}^*,\mathbf{Z}^*,\bar{\lambda}^*)$ of $\mathfrak{M}^e$ \eqref{MI} has the $\mathbf{x}^*$ component as a variational GNE of game \eqref{GM} when $X=X^e$.
\end{thm}
{\bf Proof:}
We write Algorithm \ref{alg_d} in a compact form.
Due to proximal terms $\frac{1}{2}||x_{i}-x_{i,k} ||_{R_i}^2$ and Assumption $\ref{assum1}$, the subgame in {\bf Step 1}  has a strongly monotone pseudo-subdifferential, hence its NE $\hat{\mathbf{x}}_k$ exists and is also unique.
Therefore, $\hat{x}_{i,k}=\arg\min_{x_i\in \Omega_i } \tilde{f}_i(x_i,\hat{\mathbf{x}}_{-i,k})$, and its  KKT condition
 is
\begin{equation}
\begin{array}{l} \nonumber
 \mathbf{0} \in N_{\Omega_i}(\hat{x}_{i,k})+ \partial_i  f_i(\hat{x}_{i,k},\hat{\mathbf{x}}_{-i,k})+ R_i(\hat{x}_{i,k}-x_{i,k})\\
 \quad + A^T_i[\lambda_{i,k}+H_i(A_ix_{i,k}+\sum_{l\in \mathcal{E}_i} V_{il}z_{l,k}-b_i)].
 \end{array}
\end{equation}
Concatenating all KKT conditions together and using the compact notations defined before, yields for $\hat{\mathbf{x}}_k$
\begin{equation}
\begin{array}{l}\label{equ_thm_limiting_3}
 \mathbf{0 } \in  N_{\bar{\Omega}} (\hat{\mathbf{x}}_k) + \partial  F(\hat{\mathbf{x}}_{k})+ R(\hat{\mathbf{x}}_k-\mathbf{x}_{k})\\
\qquad +\Lambda^T[\bar{\lambda}_k+H(\Lambda \mathbf{x}_{k} + \bar{V}\mathbf{Z}_k-\bar{b})].
 \end{array}
\end{equation}
We also have  $||\tilde{\mathbf{x}}_k-\hat{\mathbf{x}}_k || \leq \mu_k$ and $\mathbf{x}_{k+1}=\mathbf{x}_{k}+\rho(\tilde{\mathbf{x}}_k-\mathbf{x}_{k})$.

Let $\tilde{\lambda}_{i,k}=\lambda_{i,k}+ H_i(A_i\tilde{x}_{i,k}+\sum_{l\in \mathcal{E}_i} V_{il}z_{l,k}-b_i)$ and $\bar{\tilde{\lambda}}_k=col(\tilde{\lambda}_{1,k},\cdots,\tilde{\lambda}_{N,k})$. The compact form of {\bf Step 2} is
\begin{equation}
\begin{array}{ll}\label{equ_thm_limiting_4}
\bar{\tilde{\lambda}}_k &= \bar{\lambda}_k+H(\Lambda \tilde{\mathbf{x}}_{k} + \bar{V} \mathbf{Z}_k - \bar{b} ), \\
\bar{\lambda}_{k+1}     &= \bar{\lambda}_k+\rho(  \bar{\tilde{\lambda}}_k- \bar{\lambda}_k).
\end{array}
\end{equation}

Noticing that $\frac{1}{\rho}\lambda_{i,k+1}+\frac{\rho-1}{\rho}\lambda_{i,k} = \tilde{\lambda}_{i,k}$, we have $s_{i,k}=\tilde{\lambda}_{i,k}+H_i(A_i\tilde{x}_{i,k}+\sum_{l\in \mathcal{E}_i}V_{il}z_{l,k}-b_i)$. Denote $\bar{s}_{k}=col(s_{1,k},\cdots, s_{N,k})$, then $\bar{s}_k=\tilde{\bar{\lambda}} + H(\Lambda \tilde{\mathbf{x}}_k+\bar{V}\mathbf{Z}_{k}-\bar{b})$.
Denote $\tilde{z}_{l,k} = z_{l,k} -W_i [s_{j,k}-s_{i,k}], e_l=i \rightarrow j,$
then
$z_{l,k+1}=z_{l,k}+\rho(\tilde{z}_{l,k}-z_{l,k}).$
Denote $\tilde{\mathbf{Z}}_k=col(\tilde{z}_{1,k},\cdots,\tilde{z}_{M,k})$, then
$\tilde{\mathbf{Z}}_k  =\mathbf{Z}_k - W
\bar{V}^T\bar{s}_k $.
The updates $z_{l,k}$ are in compact form
\begin{equation}
\begin{array}{ll}\label{equ_thm_limiting_5}
\tilde{\mathbf{Z}}_k & =\mathbf{Z}_k - W\bar{V}^T [\tilde{\bar{\lambda}} + H(\Lambda \tilde{\mathbf{x}}_k+\bar{V}\mathbf{Z}_{k}-\bar{b}) ], \\
\mathbf{Z}_{k+1}     &=\mathbf{Z}_k+\rho(\tilde{\mathbf{Z}}_k-\mathbf{Z}_k).
\end{array}
\end{equation}
Using  \eqref{equ_thm_limiting_3}, \eqref{equ_thm_limiting_4}, \eqref{equ_thm_limiting_5},
Algorithm \ref{alg_d} is written compactly as 
\begin{equation}\label{alg_compact}
\begin{array}{l}
R\mathbf{x}_k -\Lambda^T [\bar{\lambda}_k+H(\Lambda \mathbf{x}_k+\bar{V}\mathbf{Z}_k-\bar{b})]
 \in (N_{\bar{\Omega}}+ \partial F+R)(\hat{\mathbf{x}}_{k}) \\
 ||\tilde{\mathbf{x}}_{k}- \hat{\mathbf{x}}_k||\leq \mu_k \\
\tilde{\bar{\lambda}}_{k}=\bar{\lambda}_k+H(\Lambda \tilde{\mathbf{x}}_{k}+\bar{V}\mathbf{Z}_k-\bar{b})\\
 \tilde{\mathbf{Z}}_{k}=\mathbf{Z}_k-W\bar{V}^T(\tilde{\bar{\lambda}}_{k}
+H(\Lambda \tilde{\mathbf{x}}_{k}+\bar{V}\mathbf{Z}_k-\bar{b})\\
\mathbf{x}_{k+1}=\mathbf{x}_k + \rho(\tilde{\mathbf{x}}_{k}-\mathbf{x}_k),\;
\bar{\lambda}_{k+1}=\bar{\lambda}_k       + \rho(\tilde{\bar{\lambda}}_{k}-\bar{\lambda}_k) \\
\mathbf{Z}_{k+1}=\mathbf{Z}_k+\rho(\tilde{\mathbf{Z}}_{k}-\mathbf{Z}_k)
\end{array}
\end{equation}
We verify next that  any limiting point of Algorithm \ref{alg_d}, or \eqref{alg_compact}, is a zero  of operator $\mathfrak{M}^e$, \eqref{MI}.
Since $\{\mu_k\}$ satisfies $\sum_{k=1}^{\infty} \mu_k < \infty$ and $\mu_k \geq 0$,  we have $\mu_k\rightarrow 0$ as $k\rightarrow \infty$.
Assume  \eqref{alg_compact} has a limiting point $col(\mathbf{x}^*,\mathbf{Z}^*,\bar{\lambda}^*)$, then we have $\mathbf{x}_{k+1}=\mathbf{x}_k=\tilde{\mathbf{x}}_k=\hat{\mathbf{x}}_k=\mathbf{x}^*$,
$\bar{\lambda}_{k+1}=\bar{\lambda}_k=\tilde{\bar{\lambda}}_k=\bar{\lambda}^*$, and
$\mathbf{Z}_{k+1}=\mathbf{Z}_k=\tilde{\mathbf{Z}}_k=\mathbf{Z}^*$.
By \eqref{alg_compact}, $col(\mathbf{x}^*,\mathbf{Z}^*,\bar{\lambda}^*)$ satisfies\vspace{-0.2cm}
\begin{align}
& -\Lambda^T [\bar{\lambda}^* + H(\Lambda \mathbf{x}^*+\bar{V}\mathbf{Z}^*-\bar{b})]
 \in (N_{\bar{\Omega}}+ \partial F)({\mathbf{x}}^*) \label{limiting_equ_1}\\
&{\bar{\lambda}^*}=\bar{\lambda}^*+H(\Lambda {\mathbf{x}^*}+\bar{V}\mathbf{Z}^*-\bar{b}) \label{limiting_equ_2}\\
&{\mathbf{Z}^*}=\mathbf{Z}^*-W\bar{V}^T({\bar{\lambda}^*}
+H(\Lambda {\mathbf{x}^*}+\bar{V}\mathbf{Z}^*-\bar{b}) \label{limiting_equ_3}
\end{align}
Since $H$, $R$ and $W$ are symmetric positive definite, \eqref{limiting_equ_2} implies that $\mathbf{0}=\Lambda {\mathbf{x}^*}+\bar{V}\mathbf{Z}^*-\bar{b}$, i.e., $A_ix_i^*-b_i = \sum_{l\in \mathcal{E}_i} V_{il}z^*_l,\forall i\in \mathcal{N}$.
Then \eqref{limiting_equ_1} and \eqref{limiting_equ_3} imply $\mathbf{0}\in \Lambda^T \bar{\lambda}^*+(N_{\bar{\Omega}}+ \partial F)({\mathbf{x}}^*)$ and
$ \mathbf{0}= \bar{V}^T{\bar{\lambda}^*}$. Using  \eqref{MI} for operator $\mathfrak{M}^e$, it follows that any limit point  of Algorithm \ref{alg_d} belongs to $zer\mathfrak{M}^e$.

 We show that any $col(\mathbf{x}^*,\mathbf{Z}^*,\bar{\lambda}^*)\in zer\mathfrak{M}^e$ has $\mathbf{x}^*$  as a variational GNE of game \eqref{GM}.
Since $\mathcal{G}$ is undirected and connected, $\bar{V}^T\bar{\lambda}^*=\mathbf{0}$ implies
$\bar{\lambda}^*= \mathbf{1}_N\otimes \lambda^* $, $\lambda^*\in \mathbf{R}^m$.
Using $\mathbf{1}^T_N V=\mathbf{0}^T_M$ and $-\Lambda \mathbf{x}^*-\bar{V} \mathbf{Z}^*+\bar{b}= \mathbf{0}$, $\mathbf{1}_N^T\otimes I_m(-\Lambda \mathbf{x}^*-V\otimes I_m \mathbf{Z}^*+\bar{b})=\mathbf{0}$ implies
$\sum_{i=1}^N A_ix^*_i =\sum_{i=1}^N b_i.$
Moreover,  $\mathbf{0} \in   \Lambda^T\bar{\lambda}^*+  (N_{\bar{\Omega}}+\partial F)(\mathbf{x}^*)$ and $\bar{\lambda}^*= \mathbf{1}_N\otimes \lambda^*$ imply
$\mathbf{0} \in A_i^T\lambda^*+ N_{\Omega_i}(x_i^*)+\partial_i f_i(x^*_i,\mathbf{x}^*_{-i}),\; \forall i\in \mathcal{N}.$
Therefore,  $\mathbf{x}^*$ and $\lambda^*$ satisfy  the KKT condition \eqref{kkt2} for the GVI \eqref{vi_original}, hence, $\mathbf{x}^*$ is a variational GNE of game  \eqref{GM} with $X=X^e$. , and all players have the same local multipliers, i.e., $\lambda_i^*=\lambda^*, \forall i\in \mathcal{N}$.

\hfill $\Box$
\subsection{Proximal parallel splitting algorithm for $X=X^i$}\label{sec_parallel_splitting}
The distributed variational GNE computation algorithm for game \eqref{GM} when $X=X^i$  is given as follows.
\begin{alg}\label{alg_d_INE}
\quad \\
\noindent\rule{0.49\textwidth}{0.7mm}
{\bf Step 1a--update of $x_{i,k}$}:
\begin{itemize}
 \item Construct a subgame where player $i$ has a decision $ x_i\in \Omega_i $ and an objective function $\tilde{f}_i(x_i,\mathbf{x}_{-i})$,\vspace{-0.2cm}
\begin{equation}\label{subgame_ine}
\tilde{f}_i(x_i,\mathbf{x}_{-i})=f_i(x_i,\mathbf{x}_{-i})+ \frac{1}{2} || x_i-x_{i,k}||_{R_i}^2
 + \lambda_{i,k}^T A_i x_i,
 \end{equation}
and denote its  NE by $\hat{\mathbf{x}}_{k}=col(\hat{x}_{1,k},\cdots, \hat{x}_{N,k})$.

\item Compute an inexact solution $\tilde{\mathbf{x}}_{k}=col(\tilde{x}_{1,k},\cdots,\tilde{x}_{N,k})$ to game in \eqref{subgame_ine} such that  $||\tilde{\mathbf{x}}_{k}-\hat{\mathbf{x}}_{k}|| \leq \mu_{k}$.

\item  Player $i$ updates its decision $x_{i,k}$ with\vspace{-0.2cm}
\begin{equation}\label{alg_d_ine_x_dynamics}
x_{i,k+1}=x_{i,k}+\rho (\tilde{x}_{i,k}-x_{i,k}).
\end{equation}
\end{itemize}

{\bf Step 1b--update of $z_{l,k}$}:
If $e_l\in \mathcal{E}_i^{out}$, then player $i$ receives $\lambda_{j,k}$, $j\in \mathcal{N}_l \setminus \{i\}$,  and updates $z_{l,k}$ with\vspace{-0.2cm}
\begin{equation}
\begin{array}{l}\label{alg_d_ine_z_dynamics}
\tilde{z}_{l,k}=z_{l,k}-W_l (\lambda_{j,k}-\lambda_{i,k}),\\
z_{l,k+1}= z_{l,k}+\rho(\tilde{z}_{l,k}-z_{l,k}).
\end{array}
\end{equation}

{\bf Step 2--update of $\lambda_{i,k}$}: player $i$ receives $\tilde{z}_{l,k},z_{l,k}$, $e_l \in \mathcal{E}_i^{in}$.\vspace{-0.15cm}
\begin{equation}\label{alg_d_ine_lambda_dynamics}
\begin{array}{l}
\tilde{{\lambda}}_{i,k} = P^{H^{-1}_i}_{\mathbf{R}^m_{+}}[\lambda_{i,k}+ H_i(A_i(2\tilde{x}_{i,k}-x_{i,k}) \\
\qquad +\sum_{l\in \mathcal{E}_i} V_{il}(2\tilde{z}_{l,k}-z_{l,k})-b_i)],\\
\lambda_{i,k+1}=\lambda_{i,k} + \rho(\tilde{{\lambda}}_{i,k}-\lambda_{i,k}).
\end{array}
\end{equation}
\noindent\rule{0.49\textwidth}{0.7mm}
All variables have the same meaning as in Algorithm \ref{alg_d}.
\end{alg}

\begin{rem}
Compared with Algorithm \ref{alg_d}, Algorithm \ref{alg_d_INE} has a different update order, that is
$\left\{
   \begin{array}{l}
     \mathbf{x}  \\
     \mathbf{Z}
   \end{array}
 \right.\rightarrow \bar{\lambda} $ rather than $ \mathbf{x} \rightarrow  \bar{\lambda} \rightarrow \mathbf{Z} $. Algorithm \ref{alg_d_INE} is called a {\it proximal parallel splitting algorithm}  since $\mathbf{x}$ and $\mathbf{Z}$ can be updated in  parallel, and only
the update of $\lambda$ utilizes the most recent information. Another difference lies in the construction of subgame at {\bf Step 1a}, i.e., \eqref{subgame_ine}  only utilizes the proximal term and Lagrangian term without considering a (linearized) quadratic penalty term. { The quadratic term in proximal ADMM is motivated by the augmented Lagrangian method for equality constrained optimization. However, the augmented Lagrangian method for inequality constrained optimization is less understood and may involve  non-differentiable terms. Hence, augmented Lagrangian methods for distributed GNE computation of inequality constrained games is beyond the scope of this paper.}
\end{rem}

Using  the same compact notations in subsection \ref{seb_sec_proximal ADMM}, such as $\bar{\lambda}$, $\mathbf{Z}$,
$R$,
$W$,
$H$,
$\bar{V}$,
$\Lambda$,
 $\Lambda^T$, and $\bar{b}$, we give next the limiting point analysis of Algorithm \ref{alg_d_INE}.
\begin{thm}\label{thm_limiting_point_ine}
Suppose that Assumption \ref{assum1}, \ref{assum2}  hold for the game \eqref{GM} when $X=X^i$. Then any limiting  point $col(\mathbf{x}^*,\mathbf{Z}^*,\bar{\lambda}^*)$ of Algorithm $\ref{alg_d_INE}$ belongs to the zeros of operator $\mathfrak{M}^i$ defined by
\begin{equation} \label{MI_ine}
{\mathfrak{M}}^i: \left(
  \begin{array}{c}
    \mathbf{x} \\
    \mathbf{Z} \\
    \bar{\lambda} \\
  \end{array}
\right) \mapsto
\left(
  \begin{array}{c}
   \Lambda^T\bar{\lambda}+  (N_{\bar{\Omega}}+\partial F)\mathbf{x} \\
    \bar{V}^T\bar{\lambda }  \\
    -\Lambda\mathbf{x} -\bar{V}\mathbf{Z}+\bar{b}+N_{\mathbf{R}^{mN}_{+}}(\bar{\lambda})    \\
  \end{array}
\right)
\end{equation}
Meanwhile, any zero $col(\mathbf{x}^*,\mathbf{Z}^*,\bar{\lambda}^*)$ of ${\mathfrak{M}}^i$, \eqref{MI_ine} has the $\mathbf{x}^*$ component as a variational GNE of game in \eqref{GM} when $X=X^i$.
\end{thm}

{\bf Proof:}  We first write Algorithm \ref{alg_d_INE} in a compact form. Since $\hat{\mathbf{x}}_k$ is an NE of subgame \eqref{subgame_ine}, $\hat{x}_{i,k}$ is an optimal solution to
$\min_{x_i\in \Omega_i} \; f_i(x_i,\hat{\mathbf{x}}_{-i,k})+ \frac{1}{2} || x_i-x_{i,k}||_{R_i}^2
 + \lambda_{i,k}^T A_i x_i.$  Under Assumption \ref{assum1} and \ref{assum2}, its optimality condition is\vspace{0.1cm}

$\mathbf{0} \in N_{\Omega_i}(\hat{x}_{i,k})+ \partial_i f_i(\hat{x}_{i,k},\hat{\mathbf{x}}_{-i,k})+ R_i(\hat{x}_{i,k}-x_{i,k})+ A_i^T \lambda_{i,k}
$. 

\noindent
 On the other hand, $\tilde{x}=P^G_{\Omega}(x)=\arg\min_{y} (\iota_{\Omega}(y)+ \frac{1}{2}||x-y ||_G^2)$ if and only if $0 \in N_{\Omega}(\tilde{x})+G(\tilde{x}-x)$.
Therefore,  the first line of \eqref{alg_d_ine_lambda_dynamics} can be written as
$\mathbf{0} \in N_{\mathbf{R}^m_{+}}(\tilde{\lambda}_{i,k})
+ H_i^{-1}\{ \tilde{\lambda}_{i,k}- \lambda_{i,k}-H_i[A_i(2\tilde{x}_{i,k}-x_{i,k})
+\sum_{l\in \mathcal{E}_i} V_{il}(2\tilde{z}_{l,k}-z_{l,k})-b_i]\}
$. 
Hence, 
for all  players we can write in compact form, \vspace{-0.15cm}
\begin{equation}\label{alg_compact_ine}
\begin{array}{l}
\mathbf{0}\in (N_{\bar{\Omega}}+ \partial F)(\hat{\mathbf{x}}_{k})+ R(\mathbf{x}_k-\hat{\mathbf{x}}_k) +\Lambda^T \bar{\lambda}_k , \\
 ||\tilde{\mathbf{x}}_{k}- \hat{\mathbf{x}}_k||\leq \mu_k,  \quad \mathbf{x}_{k+1}    = \mathbf{x}_k + \rho(\tilde{\mathbf{x}}_{k}-\mathbf{x}_k) \\
\tilde{\mathbf{Z}}_{k}=\mathbf{Z}_k-W\bar{V}^T{\bar{\lambda}}_{k}, \quad \mathbf{Z}_{k+1}    = \mathbf{Z}_k+\rho(\tilde{\mathbf{Z}}_{k}-\mathbf{Z}_k), \quad\\
\mathbf{0}\in N_{\mathbf{R}^{mN}_{+}}(\tilde{\bar{\lambda}}_{k})+H^{-1}\big \{\tilde{\bar{\lambda}}_{k}-\bar{\lambda}_k\\
-H[\Lambda (2\tilde{\mathbf{x}}_{k}-\mathbf{x})+\bar{V}(2\tilde{\mathbf{Z}}_k-\mathbf{Z}_k)-\bar{b})]\big\}\\
\bar{\lambda}_{k+1} = \bar{\lambda}_k + \rho(\tilde{\bar{\lambda}}_{k}-\bar{\lambda}_k)
\end{array}
\end{equation}

Since $R,H,W$ are positive definite,  with similar arguments in Theorem \ref{thm_limiting_point}, it can be verified that any limiting point of \eqref{alg_compact_ine} is a zero of ${\mathfrak{M}}^i$ in \eqref{MI_ine}.

Suppose $col(\mathbf{x}^*,\mathbf{Z}^*,\bar{\lambda}^*)$  is a zero of
${\mathfrak{M}}^i$ in \eqref{MI_ine}. Then with similar arguments as in Theorem \ref{thm_limiting_point},  we obtain $\bar{\lambda}^*$ to be $1_N\otimes \lambda^*, \lambda^* \in \mathbf{R}^m$. And $\mathbf{x}^*$ together with $\lambda^*$ satisfies the first line of \eqref{kkt2_ine}.
Moreover, by $\mathbf{0} \in -\Lambda\mathbf{x}^* -\bar{V}\mathbf{Z}^*+\bar{b}+N_{\mathbf{R}^{mN}_{+}}(\bar{\lambda}^*) $ and $ \bar{\lambda}^*=1_N\otimes \lambda^*$,  there exist $v_1,v_2,\cdots,v_N \in N_{\mathbf{R}^{m}_{+}}(\lambda^*)$, such that
$\mathbf{0}_{mN} =  -\Lambda \mathbf{x}^* -{V}\otimes I_m \mathbf{Z}^* + \bar{b}+ col(v_1,\cdots,v_N). $
Multiplying both sides of above equation with $\mathbf{1}^T_N\otimes I_m$ and combining with  $\mathbf{1}^T V=\mathbf{0}^T$, we have
$\mathbf{0}_{m}   = -\sum_{i=1}^N (A_ix^*_i - b_i)+\sum_{i=1}^N v_i $
We have $\sum_{i=1}^N v_i\in N_{\mathbf{R}^{m}_{+}}(\lambda^*)$ due to $v_i \in N_{\mathbf{R}^{m}_{+}}(\lambda^*)$ and $N_{\mathbf{R}^{m}_{+}}(\lambda^*)$ is a convex cone. This implies that the second line of KKT condition \eqref{kkt2_ine} is satisfied. The conclusion follows.
\hfill $\Box$

\vskip 5mm
{
\begin{rem}\label{rem_best_response_al}
Both Algorithm \ref{alg_d} and  Algorithm \ref{alg_d_INE} are {\it double-layer algorithms} since at each outer-layer iteration, players need to compute inexactly an NE of regularized subgames with a {\it given} accuracy. Since only {\it monotonicity} is assumed here, various problems could be solved with our algorithms, but the specified choice of the inner-layer algorithm should be determined according to the problem at hand. Thus, Algorithm \ref{alg_d} and  Algorithm \ref{alg_d_INE}  are ``prototype" algorithms. The inner-layer NE seeking algorithm is not specified for the following reasons.
\begin{itemize}
\item The NE seeking algorithm can and should be tailored according to the structure of the objective functions, such as  the splitting form in Remark \ref{rem_splitting}. For example, if $\partial F(\mathbf{x})=\partial L(\mathbf{x})+ \nabla_p G(\mathbf{x})$ as in Remark \ref{rem_splitting}, the possible Lipschitz continuity of $\nabla_p G(\mathbf{x})$ should be considered when choosing the NE algorithm.
\item The subgames are regularized to have {\it strongly monotone} PS/PGs due to the proximal term $\frac{1}{2} ||x_i-x_{i,k} ||_{R_i}^2$, hence efficient NE seeking algorithms available in the literature can be used, e.g., \cite{scutari,pavel3,grammatico_1,hu,baser}.
    \end{itemize}
For example, if the objective functions satisfy the assumptions in \cite{scutari},  the asynchronous distributed best-response algorithm in \cite{scutari} could be adopted for NE seeking.
Particularly, denote $\mathcal{B}(\mathbf{x}) = col(\mathcal{B}_1(\mathbf{x}_{-1}),\cdots, \mathcal{B}_N(\mathbf{x}_{-N}))$ where $\mathcal{B}_{i}(\mathbf{x}_{-i}) = \arg\min_{x_i\in \Omega_i} \tilde{f}_i(x_i,\mathbf{x}_{-i})$ given fixed $\mathbf{x}_{-i}$.
Using  Lemma 14 of \cite{scutari},  $R_i$ can be chosen such that $\mathcal{B}(\mathbf{x})$ is a  contractive map, hence the best-response algorithm enjoys a geometric convergence rate. This is even more preferable if  $\mathcal{B}(\mathbf{x})$ has a closed form. If it does not,  since $\tilde{f}_i(x_i,\mathbf{x}_{-i})$ is  strongly convex in $x_i$ given $\mathbf{x}_{-i}$ due to  $\frac{1}{2} ||x_i-x_{i,k} ||_{R_i}^2$,   $\mathcal{B}_{i}(\mathbf{x}_{-i})$
can be computed locally with the proximal gradient method that also enjoys a geometric convergence rate.

The stopping criterion for the inner layer  should be decided after the NE seeking algorithm is selected. For example, for the best-response algorithm in \cite{scutari} a termination  criterion to meet a given solution accuracy can be determined as in Remark 18 of \cite{scutari}.

We note that  a single-layer GNE seeking algorithm has been proposed in \cite{shanbhag4}, but uses diminishing step-sizes and a coordination center. Our double-layer GNE algorithm could be preferable when  there is no central node and the subgames can be easily solved.
\end{rem}

\begin{rem}
The challenges involved in  GNE seeking of game \eqref{GM} are as follows.
Firstly, the game has {\it monotone PS} without Lipschitz continuity (or the Lipschitz constant is not known prior).
Secondly, the players can only communicate peer-to-peer to coordinate to ensure
coupling constraints, even though neither $X^e$ nor $X^i$ is available to any agent.
The key idea of the proposed algorithms, i.e., Algorithm \ref{alg_d} and \ref{alg_d_INE}, is to decompose the complicated GNE seeking into sequential NE computation of regularized subgames and local coordinations. Notice that double-layer algorithms have been adopted for GNE seeking in \cite{wangjian} and \cite{pavel2}, but  only for strongly monotone games.
The proximal terms regularize the subgame such that its NE can be much easier computed. The edge variables, motivated by network flow, \cite{low}, are introduced to assist agents to reach consensus on local multipliers and to satisfy the coupling constraints.
\end{rem}
}

\section{Convergence analysis}\label{sec_convergence_analysis}

In this section, we first show that both  Algorithm \ref{alg_d} and \ref{alg_d_INE} can be derived from  a {\it preconditioned proximal point algorithm} (PPPA) for finding zeros of monotone operators. Then, based on this relationship we prove  their convergence under a sufficient choice for
the parameters $R_i$, $H_i$ and $W_l$, $\forall i\in \mathcal{N}, l\in \mathcal{E}$.

Given a  maximally monotone operator ${\mathfrak{M}}$ and a symmetric positive definite matrix $\Phi$, the {\it inexact PPPA}
with relaxation steps for finding a zero of ${\mathfrak{M}}$ is given below.
\begin{alg}\label{alg_pppa}
\quad \\
\noindent\rule{0.49\textwidth}{0.6mm}
\begin{equation}
\begin{array}{l}\label{equ_pppa}
\Phi(\varpi_{k}-\hat{\varpi}_k) \in  {\mathfrak{M}} \hat{\varpi}_k, \; \; ||\hat{\varpi}_{k}-\tilde{\varpi}_k || \leq \nu_k, \\
\varpi_{k+1}=\varpi_k+\rho(\tilde{\varpi}_k- \varpi_k ),
\end{array}
\end{equation}
where $\nu_k>0$, $\sum_{k=1}^{\infty} \nu_k <\infty$, and $\rho\in [1,2)$.
\noindent\rule{0.49\textwidth}{0.6mm}
\end{alg}

{
\begin{rem}
The proximal  algorithm for solving $\mathbf{0}\in \mathfrak{M}(x)$ (referring to Theorem 23.41 of \cite{combettes1}) is
\begin{equation}\label{alg_proximal}
\varpi_{k+1}=R_{\mathfrak{M}}\varpi_k=({\rm Id}+\mathfrak{M})^{-1} \varpi_k.
\end{equation}
which  can be equivalently written as
$\varpi_k-\varpi_{k+1} \in  \mathfrak{M} \varpi_{k+1}$.
Intuitively speaking, when $\mathfrak{M}(\varpi)$ is a linear operator $\mathfrak{M}\varpi$, each iteration of \eqref{alg_proximal} involves computing an inverse of $I+\mathfrak{M}$.
Hence, compared with \eqref{alg_proximal},
Algorithm \ref{alg_pppa} introduces a preconditioning matrix $\Phi$, considers the inexactness when evaluating the resolvent of $\mathfrak{M}$ at some specified point, and adopts an extrapolation/relexation  step.

Particularly, the preconditioning matrix $\Phi$ plays a crucial role in our algorithm design:
\begin{itemize}
\item It adds  a proximal term  $\frac{1}{2} || x_i-x_{i,k}||_{R_i}^2$ to \eqref{subgame} and \eqref{subgame_ine} that regularizes the subgames.
 \item It helps to compute the resolvent of the linear parts of ${\mathfrak{M}}^e$ and ${\mathfrak{M}}^i$ with just one step of  local communication and local computation,  without any matrix inverse.
\end{itemize}
\end{rem}
}

The next result shows the convergence of Algorithm \ref{alg_pppa}.
\begin{thm}\label{thm_pppa_convergence}
Suppose ${\mathfrak{M}}$  is  maximally monotone, and $\Phi$ is symmetric positive definite.
Suppose   $\varpi_k$ is generated by PPPA Algorithm \ref{alg_pppa} with  $\sum_{k=1}^{\infty}\nu_k <\infty ,$ $\rho \in [1,2)$. Then $\varpi_k$ converges to $\varpi^*$ and $\varpi^* \in zer{\mathfrak{M}}$. 
\end{thm}
The proof of  Theorem \ref{thm_pppa_convergence} is adapted from \cite{eckstein}, and can be found in the Appendix.

In the next two subsections, we show the convergence of Algorithm \ref{alg_d} and \ref{alg_d_INE} by relating them to PPPA Algorithm \ref{alg_pppa},  for appropriately chosen monotone operators and preconditioning matrices, and by using Theorem \ref{thm_pppa_convergence}.
\vspace{-0.3cm}
\subsection{Convergence analysis for $X=X^e$}

We introduce two auxiliary variables $\eta\in \mathbf{R}^{mN}$ and $\theta \in \mathbf{R}^{mN}$ and denote $\varpi= col(\mathbf{x},\eta,\mathbf{Z},\theta)$. Consider another operator  $\bar{\mathfrak{M}}^e$ related to $\mathfrak{M}^e$ in \eqref{MI}, defined as  $\bar{\mathfrak{M}}^e: \varpi \mapsto $
\begin{equation}\label{MI_agumented}
\left(
  \begin{array}{cccc}
  \mathbf{0} & \Lambda^T & \mathbf{0} & -\Lambda^T \\
    -\Lambda & \mathbf{0} & -\bar{V} & \mathbf{0} \\
    \mathbf{0} & \bar{V}^T & \mathbf{0} & -\bar{V}^T         \\
    \Lambda & \mathbf{0} & \bar{V} & \mathbf{0}         \\
  \end{array}
\right) \varpi
        +\left(
           \begin{array}{c}
              (N_{\bar{\Omega}}+\partial F)\mathbf{x} \\
             \bar{b} \\
             \mathbf{0} \\
             -\bar{b} \\
           \end{array}
         \right).
\end{equation}
Define a {\it preconditioning matrix} $\Phi^e$,
\begin{equation}\label{phi_preconditioned}
\Phi^e=\left(
   \begin{array}{cccc}
   R & -\Lambda^T       &  \mathbf{0}  & \Lambda^T \\
    -\Lambda & 2H^{-1}  &\bar{ V}   & \mathbf{0} \\
    \mathbf{0}        & \bar{V}^T     & W^{-1}  & \bar{V}^T         \\
    \Lambda  & \mathbf{0}       & \bar{V}   & 2H^{-1}         \\
  \end{array}
\right)
\end{equation}
where $W=\diag\{W_1,\cdots,W_M\}$, $H=\diag\{H_1,\cdots,H_N\}$.

The following result relates Algorithm \ref{alg_d}  to the PPPA Algorithm \ref{alg_pppa} for $\mathfrak{M}=\bar{\mathfrak{M}}^e$ and $\Phi=\Phi^e$.
\begin{thm}\label{thm_equvalent}
Suppose Assumption \ref{assum1} and \ref{assum2} hold.
Denote $col(\mathbf{x}_k,\mathbf{Z}_k,\bar{\lambda}_k)$,
 $\hat{\mathbf{x}}_k$
 and $col(\tilde{\mathbf{x}}_k,\tilde{\mathbf{Z}}_k,\tilde{\bar{\lambda}}_k)$
as points generated by Algorithm \ref{alg_d}
for initial points $\mathbf{x}_0,\mathbf{Z}_0,\bar{\lambda}_0$.
Denote
$\varpi_k=col(\mathbf{x}^{'}_k,\eta_k,\mathbf{Z}^{'}_k,\theta_k)$, $\hat{\varpi}_k=col(\hat{\mathbf{x}}^{'}_k,\hat{\eta}_k, \hat{\mathbf{Z}}^{'}_k,\hat{\theta}_k)$,
and
$\tilde{\varpi}_k=col(\tilde{\mathbf{x}}^{'}_k,\tilde{\eta}_k, \tilde{\mathbf{Z}}^{'}_k,\tilde{\theta}_k)$
as the points generated by the PPPA Algorithm \ref{alg_pppa} with $\mathfrak{M}=\bar{\mathfrak{M}}^e$ and $\Phi=\Phi^e$
for initial points $\mathbf{x}^{'}_{0}=\mathbf{x}_{0}, \eta_0 = \bar{\lambda}_0+H(\Lambda \mathbf{x}_{0}+\bar{V}^T \bar{\lambda}_0-\bar{b} ),\mathbf{Z}^{'}_0=\mathbf{Z}_0,\theta_0=\mathbf{0}$.
Then,  any sequence $col(\mathbf{x}_k,\mathbf{Z}_k,\bar{\lambda}_k)$ can be derived from some sequence $\varpi_k=col(\mathbf{x}^{'}_k,\eta_k, \mathbf{Z}^{'}_k,\theta_k)$ as follows
\begin{equation}
\begin{array}{l}\label{correspondence}
\mathbf{x}_k    = \mathbf{x}^{'}_k,  \qquad \mathbf{Z}_k    = \mathbf{Z}^{'}_k,  \\
\bar{\lambda}_k = \eta_k-\theta_k-H(\Lambda \mathbf{x}^{'}_k+\bar{V}\mathbf{Z}^{'}_k-\bar{b}),
\end{array}
\end{equation}
for some nonnegative sequence $\{\nu_k\}$ such that $\sum_{k=1}^{\infty}{\nu_k}< \infty$.
\end{thm}
The proof of Theorem \ref{thm_equvalent} is based on an induction argument and is given in the Appendix.

\begin{rem}
The standard ADMM for optimization can be derived from the {\it Douglas-Rachford} (DS) splitting method for dual optimization problems, and analyzed as a proximal-point algorithm, see \cite{eckstein} and \cite{combettes1}. For proximal ADMM,  the analysis in \cite{hebingsheng} shows that  the posterior second coordinate is not available when updating the first one.
 That is the reason why we split $\bar{\lambda}$ into $\eta$ and $\theta$, to have a higher order dynamics. The {\it preconditioned DS splitting method}
 recently introduced in \cite{sunhongpeng}, might  lead to proximal ADMM.
Compared to \cite{sunhongpeng}, our algorithm applies relaxation steps to all coordinates and considers inexactness in solving the subproblems.
\end{rem}

We prove the convergence of Algorithm \ref{alg_d}, 
by exploiting the relationship given in Theorem \ref{thm_equvalent} and using  Theorem \ref{thm_pppa_convergence}.

\begin{thm}\label{thm_convergence}
Suppose Assumption \ref{assum1} and \ref{assum2} hold for  game \eqref{GM} when $X=X^e$,  and
parameters (step-sizes)  $R_i,H_i,W_l$ are  symmetric positive definite,  chosen such that    $R-\Lambda^T H\Lambda$ and $W^{-1}-\bar{V}^TH\bar{V}$ are positive definite.
Then, any $col(\mathbf{x}_k, \mathbf{Z}_k,\bar{\lambda}_k)$ generated by Algorithm \ref{alg_d} converges to $col(\mathbf{x}^*, \mathbf{Z}^*, \bar{\lambda}^*) \in zer \mathfrak{M}^e$. Furthermore, $\mathbf{x}^*$ is a variational GNE of game in \eqref{GM} when $X=X^e$, and $\bar{\lambda}^*=\mathbf{1}_N\otimes \lambda^*$, $\lambda^* \in \mathbf{R}^m$.
\end{thm}
{\bf Proof:} By Theorem \ref{thm_equvalent}, Algorithm 1 is related to PPPA Algorithm 3 for $\bar{\mathfrak{M}}^e$,  \eqref{MI_agumented}, $\Phi^e$,  \eqref{phi_preconditioned}. Convergence follows by Theorem \ref{thm_pppa_convergence} if we show that $\bar{\mathfrak{M}}^e$  is maximally monotone and   $\Phi^e$  is positive definite. 
Denote  $\varpi=col(\mathbf{x},\eta,\mathbf{Z},\theta)$, then
\begin{equation}
\begin{array}{l}\nonumber
\varpi^T \Phi^e \varpi= \mathbf{x}^T R\mathbf{x} - 2\mathbf{x}^T\Lambda^T \eta + 2\mathbf{x}^T\Lambda^T \theta
+2\eta H^{-1}\eta \\
+ 2 \eta^T \bar{V}\mathbf{Z} + \mathbf{Z}^TW^{-1}\mathbf{Z}+ 2\mathbf{Z}^T \bar{V}^T\theta+2 \theta^T H^{-1}\theta\\
= ||H\Lambda \mathbf{x}+\theta-\eta||_{H^{-1}}^2 + ||H\bar{V}\mathbf{Z}+\theta+\eta ||_{H^{-1}}^2\\
+ ||\mathbf{x} ||^2_{R-\Lambda^T H\Lambda}+ || \mathbf{Z}||^2_{W^{-1}-\bar{V}^TH\bar{V}}
\end{array}
\end{equation}
Since $R-\Lambda^T H\Lambda$ and $W^{-1}-\bar{V}^TH\bar{V}$ are positive definite, it follows immediately that 
$\Phi^e$ is positive definite.

Operator $\bar{\mathfrak{M}}^e$, \eqref{MI_agumented},  is written as the sum of two operators.
The first  is a skew-symmetric linear operator, hence, is maximally monotone with domain of whole space. $N_{\bar{\Omega}}$ is maximally monotone as a normal cone operator of a closed convex set, and $\partial F(\mathbf{x})$ is also maximally monotone by Assumption \ref{assum1}. Since their domains coincide, $N_{\bar{\Omega}}+\partial F$ is maximally monotone, and the 2$^{nd}$ term in \eqref{MI_agumented} is  maximally monotone  as the  Cartesian product of
maximally monotone operators.

By Theorem \ref{thm_equvalent},
for any sequence $col(\mathbf{x}_k,\mathbf{Z}_k,\bar{\lambda}_k)$  generated from Algorithm \ref{alg_d}, we can find $\varpi_k=col(\mathbf{x}^{'}_k,\eta_k,\mathbf{Z}^{'}_k, \theta_k)$ generated from Algorithm \ref{alg_pppa} such that \eqref{correspondence} holds for all $k$ and $\sum_{k=1}^{\infty}\nu_k<\infty$.
By Theorem \ref{thm_pppa_convergence}, $\varpi_k$ converges to $\varpi^*=col(\mathbf{x}^{'*},\eta^*,\mathbf{Z}^{'*}, \theta^*)$ and $\varpi^*\in zer\bar{\mathfrak{M}}^e$.
By \eqref{correspondence}, $col(\mathbf{x}_k,\mathbf{Z}_k,\bar{\lambda}_k)$  also converges to $col(\mathbf{x}^*,\mathbf{Z}^*,\bar{\lambda}^*)$ such that $\mathbf{x}^*=\mathbf{x}^{'*}$, $\mathbf{Z}^*=\mathbf{Z}^{'*}$, and $\bar{\lambda}^*=\eta^{*}-\theta^{*}-H(\Lambda \mathbf{x}^{'*}+\bar{V}\mathbf{Z}^{'*}-\bar{b})$.

Since $\varpi^*=col(\mathbf{x}^{'*},\eta^*,\mathbf{Z}^{'*}, \theta^*)$ is a zero of $\bar{\mathfrak{M}}^e$, \eqref{MI_agumented},
 we have $ \Lambda \mathbf{x}^{'*}+\bar{V}\mathbf{Z}^{'*} - \bar{b}=\mathbf{0}$, so that  $\bar{\lambda}^*=\eta^*-\theta^*$, and \vspace{-0.1cm}
\begin{equation}\nonumber
\mathbf{0}  \in \Lambda^T(\eta^*-\theta^*)+ (N_{\bar{\Omega}}+\partial F)\mathbf{x}^{'*}, \qquad
\mathbf{0 } \in \bar{V}^T(\eta^*-\theta^*).
\end{equation}
Using  $\mathbf{x}^{'*}=\mathbf{x}^*$, $\mathbf{Z}^{'*}=\mathbf{Z}^*$ and the definition of $\mathfrak{M}^e$ in \eqref{MI},  it follows that $col(\mathbf{x}_k,\mathbf{Z}_k,\bar{\lambda}_k)$  generated from Algorithm \ref{alg_d} converges to
$col(\mathbf{x}^*,\mathbf{Z}^*,\bar{\lambda}^*) \in zer\mathfrak{M}^e$.
By Theorem \ref{thm_limiting_point}, 
 $\mathbf{x}_k$ converges to $\mathbf{x}^*$, a variational GNE of the game in \eqref{GM}, and players'  local multipliers  converge to the same $\lambda^*$, which together with $\mathbf{x}^*$ satisfies KKT condition in \eqref{kkt2}.
\hfill $\Box$

\begin{rem}
If  $H_i$ is chosen to be a diagonal positive matrix,  $R_i $ and $W_l$ can be chosen using diagonally dominance to ensure $R-\Lambda^T H\Lambda$ and $W^{-1}-\bar{V}^TH\bar{V}$ are positive definite. In this case, the parameters $R_i,H_i$ and $W_l$ can be chosen independently by player $i$ with just local data and computation.
\end{rem}

\subsection{Convergence analysis for $X=X^i$}\label{sec_ine}

The next result shows the convergence of Algorithm \ref{alg_d_INE}. 
\begin{thm}
Suppose Assumption \ref{assum1} and \ref{assum2} hold for game \eqref{GM} when $X=X^i$, and
parameters $R_i,H_i,W_l$ are  symmetric positive definite,  such that the matrix $\Phi^i$ is positive definite,
\begin{equation}
\Phi^i=\left(
  \begin{array}{ccc}
    R & 0 & -\Lambda^T \\
    0 & W^{-1} & -V^T \\
    -\Lambda & -V & H^{-1} \\
  \end{array}
\right)
\end{equation}
Then any $col(\mathbf{x}_k, \mathbf{Z}_k,\bar{\lambda}_k)$ generated by Algorithm \ref{alg_d_INE} converges to $col(\mathbf{x}^*, \mathbf{Z}^*, \bar{\lambda}^*) \in zer {\mathfrak{M}}^i$ in \eqref{MI_ine}. Furthermore, $\mathbf{x}^*$ is a variational GNE of game in \eqref{GM} and $\bar{\lambda}^*=\mathbf{1}_N\otimes \lambda^*$, $\lambda^* \in \mathbf{R}_{+}^m$.
\end{thm}

{\bf Proof:}
Consider the PPPA Algorithm  \ref{alg_pppa} with $\varpi=col(\mathbf{x},\mathbf{Z},\bar{\lambda})$, for $\Phi=\Phi^i$ and ${\mathfrak{M}}={\mathfrak{M}}^i$ and $\nu_k=\mu_k$. After manipulations, the PPPA algorithm  gives  \eqref{alg_compact_ine}. Hence, Algorithm \ref{alg_d_INE} can be derived from Algorithm \ref{alg_pppa} via a one-to-one correspondence relation. Notice that  ${\mathfrak{M}}^i$ in \eqref{MI_ine} can be written as the sum of a skew-symmetric linear operator and a product of
$(N_{\bar{\Omega}}+\partial F)\mathbf{x}\times\mathbf{ 0} \times  N_{\mathbf{R}^{mN}_{+}}(\bar{\lambda})$.
Under  Assumption \ref{assum1} and \ref{assum2}, with similar arguments as in Theorem \ref{thm_convergence},
we can show that  ${\mathfrak{M}}^i$ in \eqref{MI_ine} is maximally monotone.
Since $\Phi^i$ is symmetric positive definite, 
by Theorem \ref{thm_pppa_convergence}, PPPA Algorithm \ref{alg_pppa} converges. Therefore, Algorithm \ref{alg_d_INE} converges to a zero of $\mathfrak{M}^i$, 
and the conclusion follows by invoking Theorem \ref{thm_limiting_point_ine}.
\hfill $\Box$

\begin{rem}
Our recent work \cite{yipeng}  considers  GNE computation for games with inequality affine constraints, but assumes  a strongly monotone and Lipschitz continuous PG, with inertial steps for possible acceleration.
In this paper, we only assume a monotone PS, consider the inexactness when solving subproblems, and use  relaxation steps  for possible acceleration.
{  Moreover, as seen in the convergence analysis,
 both Algorithm \ref{alg_d} and \ref{alg_d_INE}  can be regarded as fixed-point iterations for averaged operators, hence the convergence rate for fixed-point residuals could be derived based on an  analysis as in \cite{wotao}.}
\end{rem}

\section{Application and simulation studies}\label{sec_simulations}

\subsection{Rate control game over wireless ad-hoc networks}

\begin{figure}
  \centering
  \includegraphics[width=3in]{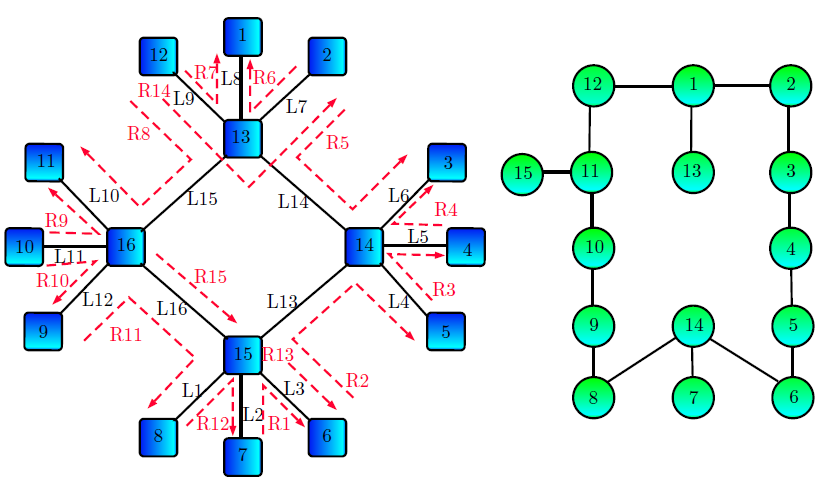}\\
  \caption{(a): Wireless Ad-Hoc Network. (b): Communication graph.\label{fig_wanet}}
\end{figure}

This example is adapted from  \cite{pavel4}.
Consider a wireless ad-hoc network (WANET)
with $16$ nodes and $16$  links $\{L_1,\cdots, L_{16}\}$ as shown in  Fig. \ref{fig_wanet}. 
There are $15$ users $\{U_1,...,U_{15}\}$ who want to transfer data through the links.
$R_i$ is the path adopted by user $U_i$, and $L_j\in R_i$ if user $U_i$ transfers data through link $L_j$.
User $U_i$ decides its data rate $x_i$, and { should satisfy a local constraint $0\leq x_i\leq B_i$}.
In Fig. \ref{fig_wanet}, the solid lines represent the links
$\{L_1,\cdots, L_{16}\}$, and dashed line displays each path $R_i$.
Denote $A=[A_1,\cdots, A_{15}]\in \mathbf{R}^{16\times15}$ where $A_i\in \mathbf{R}^{16}$, and $A_i$ has its $j$th element to be $1$ if $U_i$ uses $L_j$ and  to be $0$, otherwise.
Link $L_j$ has a maximal capacity $C_j>0$.
Denote $C=col(C_1,\cdots,C_{16})$, hence all users' data rate $\mathbf{x}$ should satisfy the inequality coupling constraint  $A\mathbf{x}\leq C$.
The objective function of user $U_i$ is
$f_i(x_i,\mathbf{x}_{-i}) = -u_i(x_i) + D^T(\mathbf{x}) A_i x_i
$,
where { $u_i(x_i)=\chi_i \log(x_i+1)$ is user $i$'s utility function}, and
$D(\mathbf{x})=col(d_1(\mathbf{x}),\cdots, d_{16}(\mathbf{x}))$ with $d_j(\mathbf{x}) = \frac{\kappa_j}{C_j-[A\mathbf{x}]_j + \xi_j}$ maps $\mathbf{x}$ to the unit delays of each link.
The parameters are randomly drawn as follows:
$C_j\in [10,15]$, $B_i \in [5,10]$, $\chi_i \in [10,20]$, $\kappa_j \in [10,30]$ and $\xi_j \in [20,40]$, { and are numerically verified to ensure Assumption \ref{assum1}.}

\begin{figure}
  \centering
  \includegraphics[width=3.5in]{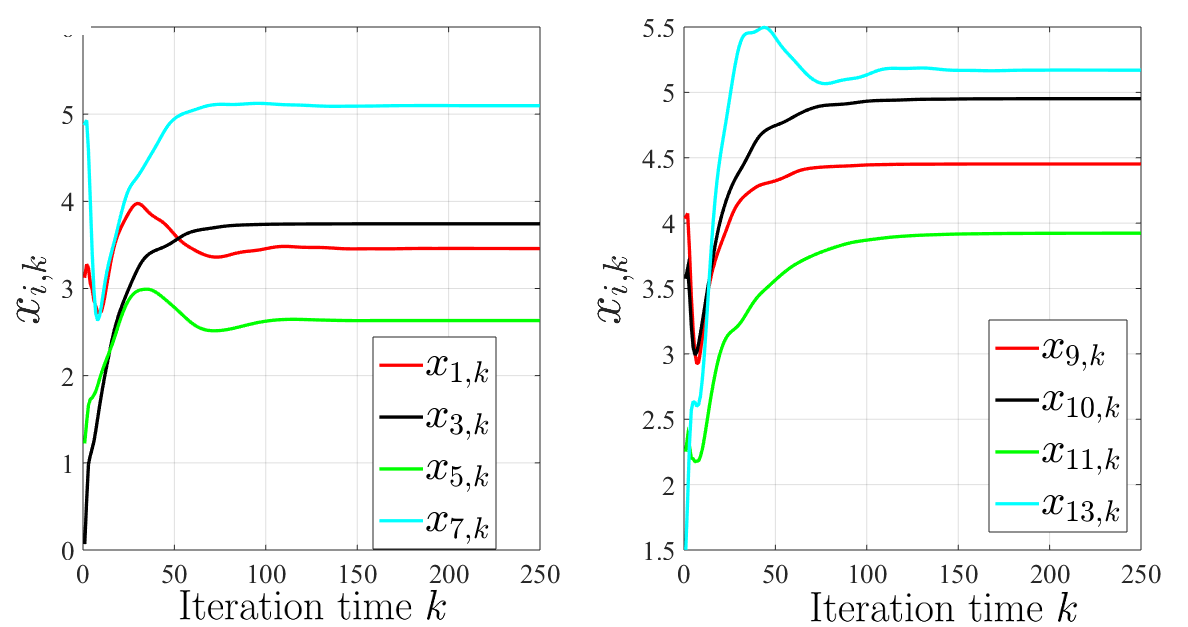}\\
  \caption{  The trajectories of selected users' data rate $x_{i,k}$, which  show the convergence of Algorithm \ref{alg_d_INE}.\label{fig_sim_1}  }
\end{figure}

\begin{figure}
  \centering
  \includegraphics[width=3.5in]{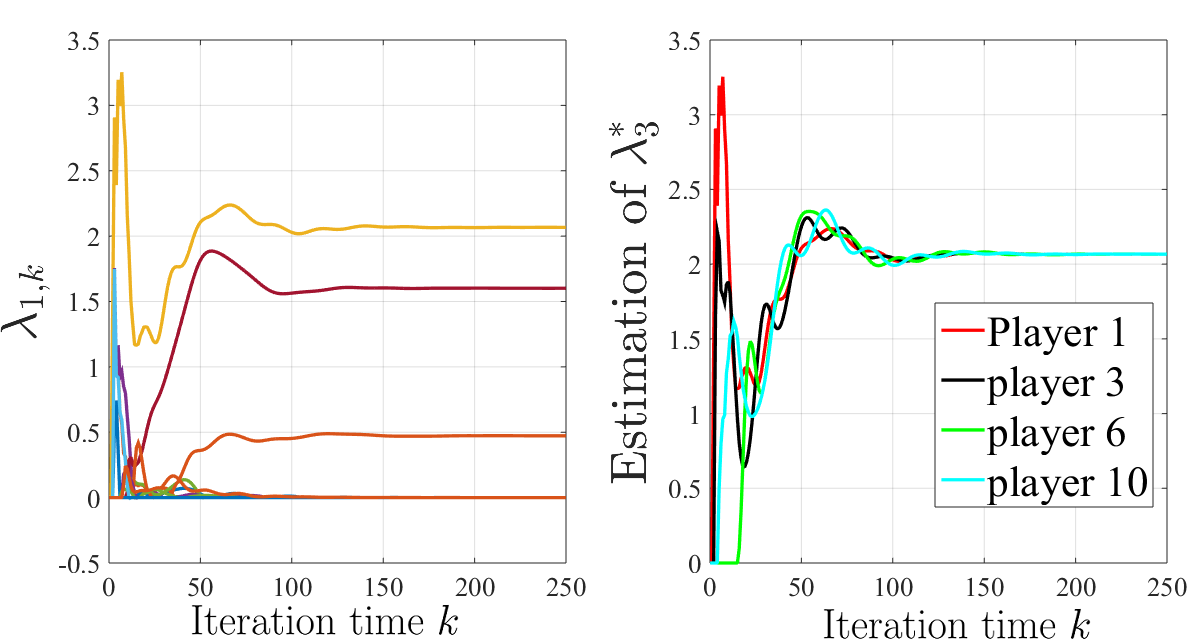}\\
  \caption{ (a): The trajectories of local multiplier $\lambda_{1,k}$ of player $1$. (b):  The trajectories of selected users' estimations of the third component of $\lambda^*$. It shows that all the players find the same multiplier $\lambda^*$. \label{fig_sim_2}  }
\end{figure}

\begin{figure}
  \centering
  \includegraphics[width=3.5in]{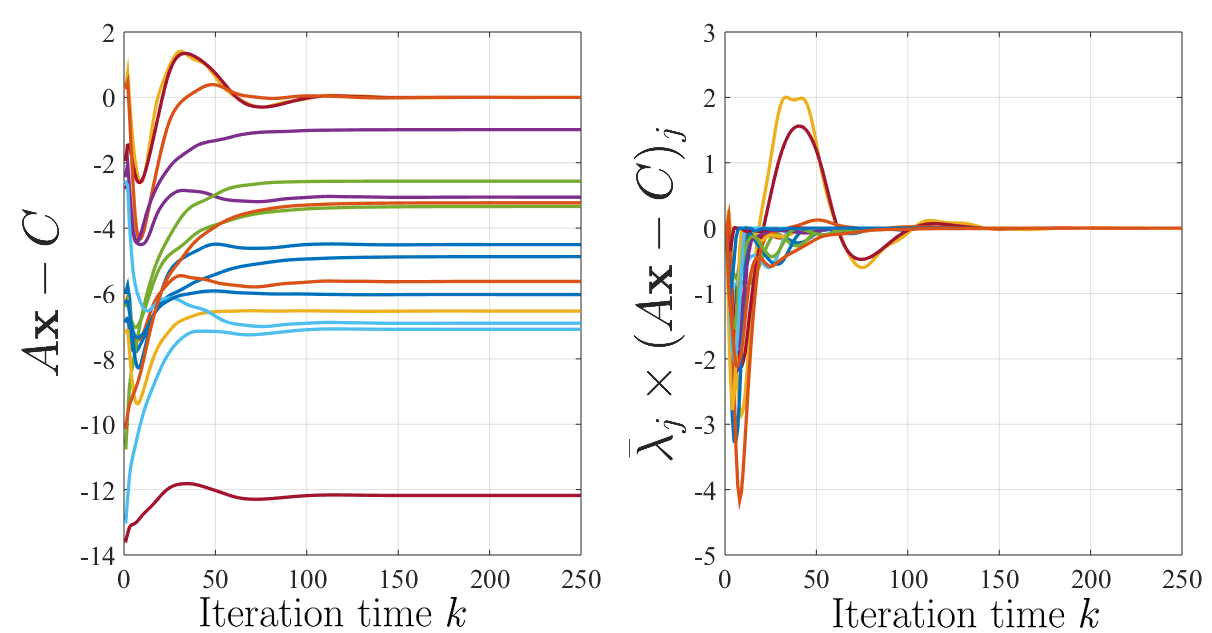}\\
  \caption{ (a) The trajectories of violations of the coupling constraint $A\mathbf{x}-C$.
  (b)$\bar{\lambda}_j$ is the averaging of the $j$th component of all players' local  multipliers.
  It shows that the coupling constraint is asymptotically satisfied, and the complementary condition $ \lambda^* \bot A\mathbf{x}^*-C$  is asymptotically satisfied. \label{fig_sim_3}  }
\end{figure}

 We use Algorithm \ref{alg_d_INE}. Each player has a local $C_i=\frac{1}{15}C$, and  has local step-sizes
 $R_i = 10$, $H_i = 0.5I_{16}$, $W_l = 0.5I_{16}$ and $\rho=1.1$.
Players communicate over the graph in Fig. \ref{fig_wanet}, with edges arbitrarily ordered.
The initial point ${x}_{i,0}$ is randomly chosen within $[0,B_i]$, and initial $\lambda_i$, $z_l$ are chosen to be zero.
The subgames are solved using gradient methods in \cite{pavel3} to get the exact NE $\hat{\mathbf{x}}_k$, and each $\tilde{\mathbf{x}}_k$ is chosen as the first point  satisfying $||\tilde{\mathbf{x}}_k- \hat{\mathbf{x}}_k|| < \frac{1}{k^2}$. The simulation results are shown in Fig. \ref{fig_sim_1}-\ref{fig_sim_3}.

\vspace{-0.3cm}

\subsection{Task allocation game}

\begin{figure}
\begin{center}
\begin{tikzpicture}[->,>=stealth',shorten >=0.3pt,auto,node distance=2cm,thick,
  rect node/.style={rectangle, ball color={rgb:red,0;green,0.2;yellow,1},font=\sffamily,inner sep=1pt,outer sep=0pt,minimum size=14pt},
  wave/.style={decorate,decoration={snake,post length=0.1mm,amplitude=0.5mm,segment length=3mm},thick},
  main node/.style={shape=circle, ball color=green!20,text=black,inner sep=1pt,outer sep=0pt,minimum size=15pt},scale=0.88]


  \foreach \place/\i in {{
  (-3,0.5)/1},
  {(-3,-0.5)/2},
  {(-3,-1.5)/3},
  {(-1,-0.5)/4},
  {(-1,-1.5)/5},
  {(1,0.5)/6},
  {(1,-0.5)/7},
  {(3,0.5)/8},
  {(3,-0.5)/9},
  {(3,-1.5)/10},
  {(0,1)/11},
  {(0,0)/12},
  {(0,-1)/13},
  {(0,-2)/14}}
    \node[main node] (a\i) at \place {};

      \node at (-3,0.5) {\rm \color{black}{$w_{1}$}};
      \node at (-3,-0.5){\rm \color{black}{$w_2$}};
      \node at (-3,-1.5){\rm \color{black}{$w_3$}};
      \node at (-1,-0.5){\rm \color{black}{$w_{4}$}};
      \node at (-1,-1.5){\rm \color{black}{$w_{5}$}};
      \node at (1,0.5)  {\rm \color{black}{$w_{6}$}};
      \node at (1,-0.5) {\rm \color{black}{$w_7$}};
      \node at (3,0.5)  {\rm \color{black}{$w_8$}};
      \node at (3,-0.5) {\rm \color{black}{$w_9$}};
      \node at (3,-1.5) {\rm \color{black}{$w_{10}$}};
      \node at(0,1)     {\rm \color{black}{$w_{11}$}};
      \node at (0, 0)   {\rm \color{black}{$w_{12}$}};
      \node at (0,-1)   {\rm \color{black}{$w_{13}$}};
      \node at (0,-2)   {\rm \color{black}{$w_{14}$}};

  \foreach \place/\x in {{(-2,1)/1},{(-2,0)/2},{(-2,-1)/3},
    {(-2,-2)/4}, {(2,1)/5}, {(2,0)/6}, {(2,-1)/7}, {(2,-2)/8}}
  \node[rect node] (b\x) at \place {};

      \node at (-2,1) {\rm \color{black}{$T_1$}};
      \node at (-2,0){\rm \color{black}{$T_2$}};
      \node at (-2,-1)  {\rm \color{black}{$T_3$}};
      \node at (-2,-2) {\rm \color{black}{$T_4$}};
      \node at (2,1){\rm \color{black}{$T_5$}};
      \node at (2,0) {\rm \color{black}{$T_6$}};
      \node at (2,-1)  {\rm \color{black}{$T_7$}};
      \node at (2,-2) {\rm \color{black}{$T_8$}};

  \path[dashed,->,blue,thick]               (a1) edge (b1);
  \path[->,red,thick]                       (a1) edge (b2);

  \path[dashed,->,blue,thick]               (a2) edge (b2);
  \path[->,red,thick]                       (a2) edge (b3);

  \path[dashed,->,blue,thick]               (a3) edge (b3);
  \path[->,red,thick]                       (a3) edge (b4);

  \path[dashed,->,blue,thick]               (a4) edge (b2);
  \path[->,red,thick]                       (a4) edge (b3);

  \path[dashed,->,blue,thick]               (a5) edge (b3);
  \path[->,red,thick]                       (a5) edge (b4);

  \path[dashed,->,blue,thick]               (a6) edge (b5);
  \path[->,red,thick]                       (a6) edge (b6);

  \path[dashed,->,blue,thick]               (a7) edge (b6);
  \path[->,red,thick]                       (a7) edge (b7);

  \path[dashed,->,blue,thick]               (a8) edge (b5);
  \path[->,red,thick]                       (a8) edge (b6);

  \path[dashed,->,blue,thick]               (a9) edge (b6);
  \path[->,red,thick]                       (a9) edge (b7);

  \path[dashed,->,blue,thick]               (a10) edge (b7);
  \path[->,red,thick]                       (a10) edge (b8);

  \path[dashed,->,blue,thick]               (a11) edge (b2);
  \path[->,red,thick]                       (a11) edge (b5);

  \path[dashed,->,blue,thick]               (a12) edge (b2);
  \path[->,red,thick]                       (a12) edge (b6);

  \path[dashed,->,blue,thick]               (a13) edge (b3);
  \path[->,red,thick]                       (a13) edge (b7);

  \path[dashed,->,blue,thick]               (a14) edge (b4);
  \path[->,red,thick]                       (a14) edge (b7);
\end{tikzpicture}
\end{center}
\caption{Task allocation game: An edge from $w_i$ to $T_j$ on this graph implies that a part of  worker $w_i$'s output is  allocated to task $T_j$.}\label{fig_task_allocation_game}
\end{figure}
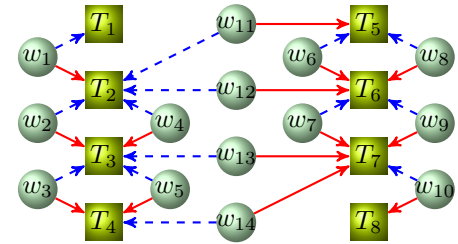

In this part, we consider a task allocation game with $8$ tasks $\{T_1,\cdots,T_8\}$ and $14$ processors (workers) $\{w_1,\cdots,w_{14}\}$.
Each task $T_j$ is quantified as a load of $C_j>0$ that should be met by the workers.
{ Each worker $w_i$ decides its working output $x_i=col(x^1_{i},x^2_{i},x^3_{i}, x^4_{i})\in \mathbf{R}^4$ within its capacity $\mathbf{0} \leq x_i\leq B_i, B_i \in \mathbf{R}_{+}^4$.
If worker $w_i$ allocates a part of its output to task $T_j$, there is an arrow $w_i\rightarrow T_j$ in Fig. \ref{fig_task_allocation_game}, either blue or red.
Specifically, if $w_i$ allocates $x^1_{i},x^2_{i}$ to $T_j$, there is a dashed blue arrow in Fig. \ref{fig_task_allocation_game}, and
if $w_i$ allocates $x^3_{i},x^4_{i}$ to $T_j$, there is a solid red arrow in Fig. \ref{fig_task_allocation_game}
Define a matrix $A=[A_1, \cdots, A_{15}]\in \mathbf{R}^{8\times 56}$
with $A_i=[a^1_i,a_i^2,a_i^3,a_i^4] \in \mathbf{R}^{8\times 4}$ quantifying how the output of worker $w_i$  is allocated to each task. Each column $a_i^k$ has only one element being nonzero, and the $j$th element of $a_i^1$ or $a_i^2$ is nonzero if there is a  dashed blue arrow $w_i\rightarrow T_j$ on Fig. \ref{fig_task_allocation_game},
and the $j$th element of $a_i^3$ or $a_1^4$ is nonzero if there is a  red arrow $w_i\rightarrow T_j$ on Fig. \ref{fig_task_allocation_game}.
The nonzero elements in $A_i$ are randomly chosen from $[0.5,1]$.
It is required that
the tasks should be met by the working output of the players.
Denote $C=col(C_1,\cdots,C_8)$, then the workers have an equality coupling constraint:
$A\mathbf{x}=C$.
The objective function of player (worker) $w_i$ is
$f_i(x_i,\mathbf{x}_{-i}) = c_i(x_i) - R^T(\mathbf{x}) A_i x_i$.
Here,  $c_i(x_i)$ is a cost function of worker $w_i$ and is taken as
$c_i(x_i)=\sum_{s=1}^4 \max\{ q^s_i {x^s_i}^2-\xi^s_i x^s_i, l^s_i x^s_i\} +  (p_i^Tx_i-d_i)^2+x_i^T S_ix_i$.
$R(\mathbf{x})=col(R_1(\mathbf{x}),\cdots,R_8(\mathbf{x}))$
is a vector function that maps the workers' output  to the award price of each task, and
$R_j(\mathbf{x})=\kappa_j - \chi_j \log( [A\mathbf{x}]_j + 1 )$.
Parameters of the problem are randomly drawn as follows:
$C_j\in [1,2]$, $\chi_j \in [0.1,0.6]$, $\kappa_j \in [10,20]$, $q^s_i\in [1,2]$, $\xi^s_i \in [6,12]$, $d_i\in [1,2]$,  and $l^s_i \in [1,3]$.
$p_i \in \mathbf{R}^{4}$ is a randomly generated stochastic vector,  $S_i\in \mathbf{R}^{4\times 4}$ is a randomly generated positive definite matrix,  and each element of $B_i$ is  drawn from $[1,3]$.
The parameters are numerically checked to ensure Assumption \ref{assum1}.

We apply  Algorithm \ref{alg_d} to this problem, over a communication graph as in Fig. \ref{fig_wanet}, without node 15 and its adjacent edge, and with the remaining edges  arbitrarily ordered.
Each player has a local $C_i=\frac{1}{15}C$, and local step-sizes
 $R_i$, $H_i$, $W_l$ that are all diagonal matrices with nonzero elements uniformly drawn from $[4, 8]$, $[0.2, 0.4]$ and [0.2, 0.4] respectively. The relaxation step-size is taken as $\rho=1.1$.
The initial ${x}_{i,0}$ is randomly chosen within $0\leq x_{i,0}\leq B_i$, and initial $\lambda_i$, $z_l$ are chosen to be zeros.
The subgames are solved using subgradient methods in \cite{shanbhag} to get the exact NE $\hat{\mathbf{x}}_k$, and each $\tilde{\mathbf{x}}_k$ is chosen to be the first point on the trajectory satisfying $||\tilde{\mathbf{x}}_k- \hat{\mathbf{x}}_k|| < \frac{1}{k^2}$. The simulation results are shown in Fig. \ref{fig_sim_4}-\ref{fig_sim_5}.

\begin{figure}
  \centering
  \includegraphics[width=3.3in]{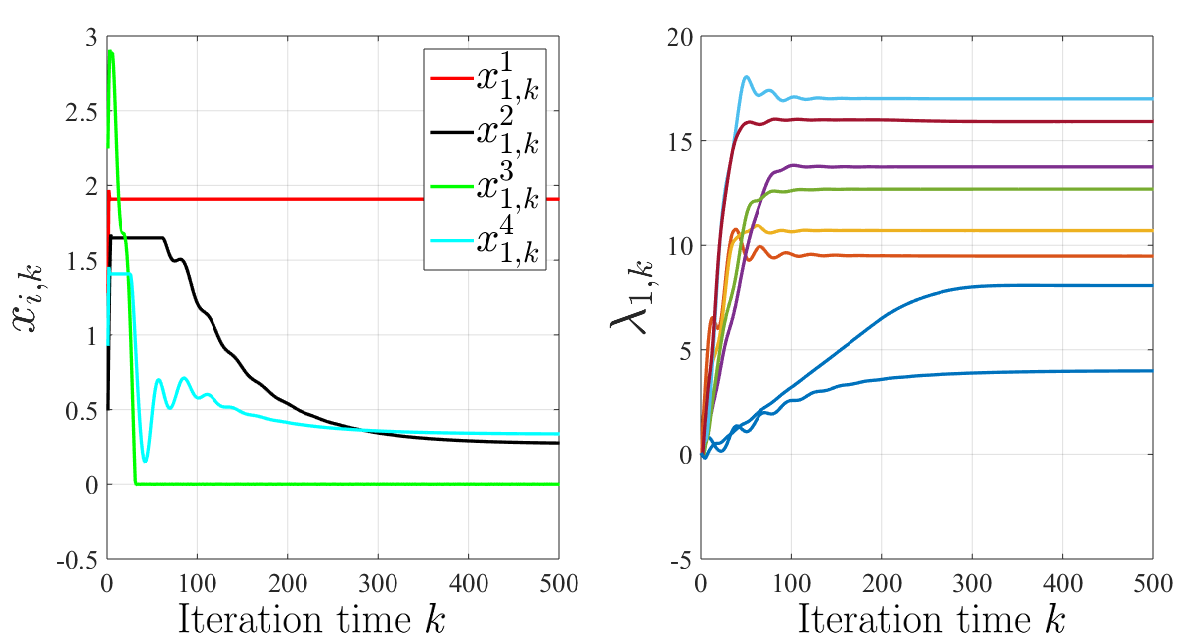}\\
  \caption{ (a): The trajectories of player $1$'s working output $x^j_{1,k},j=1,\cdots,4$, which  show the convergence of Algorithm \ref{alg_d}. (b): The trajectories of player $1$'s multiplier $\lambda_{1,k}$.\label{fig_sim_4}  }
\end{figure}
\begin{figure}
  \centering
  \includegraphics[width=3.3in]{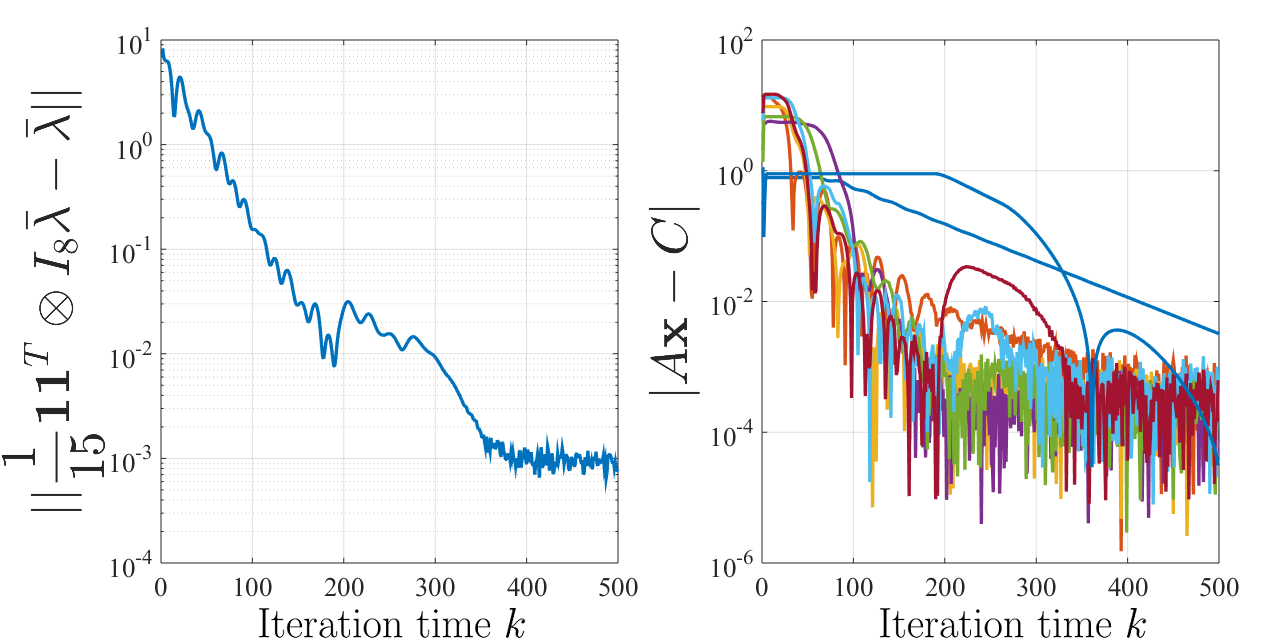}\\
  \caption{  (a):  The trajectories of the consensus errors of local multipliers.  (b) The trajectories of violations of the coupling constraint $A\mathbf{x}=C$\label{fig_sim_5}  }
\end{figure}

}
\section{Conclusions}\label{sec_concluding}

In this paper, we considered  GNE computation of monotone games with affine coupling constraints.
We proposed center-free distributed algorithms for both equality and inequality constraints, based on a preconditioned
proximal point algorithm. We decomposed the GNE computation  into sequential NE computation of regularized subgames and local coordination of multipliers and auxiliary variables. We considered inexactness in solving the subgames   and  incorporated relaxation steps.
We proved their convergence by resorting to the theory of proximal algorithms and averaged operators.

There are still a lot of promising open problems.
Motivated by \cite{pavel3} and \cite{pavel4}, it is appealing to consider distributed GNE seeking when players cannot observe all other players' decisions.
Motivated by \cite{liangshu,grammatico_2} and \cite{lygeros2}, center-free GNE seeking of monotone aggregative games with discrete-time algorithms is still open.
It is appealing to develop asynchronous distributed GNE computation algorithms with delayed information,
and consider the problem when the players interact over switching and directed communication graphs. { As important is to consider computational GNE seeking algorithms together with the mechanism design which can ensure that players faithfully report their states and auxiliary variables, possibly by providing proper incentive or punishment.}

\section*{Appendix}

Essentially, the proof of Theorem \ref{thm_pppa_convergence} utilizes the following facts: $\Phi^{-1}{\mathfrak{M}}$ is  maximally monotone under the $\Phi-$induced norm $||\cdot||_{\Phi}$;
$R_{\Phi^{-1}\mathfrak{M}}$ is a $\frac{1}{2}$-averaged operator; Proposition 4.25 of \cite{combettes1} for averaged operators and Robbins-Siegmund lemma for sequence convergence, given as follows.

\begin{lem}[Proposition 23.7 of \cite{combettes1}]\label{lem_maxiamlly_monotone}
If operator $\mathfrak{A}$ is maximally monotone, then $T=R_{\mathfrak{A}}=({\rm Id}+\mathfrak{A})^{-1}$ is firmly nonexpasive, and $domR_{\mathfrak{A}}=\mathbf{R}^m$.
\end{lem}

\begin{lem} [Proposition 4.25 of \cite{combettes1}]\label{lem_alpha_opertor}
Given an operator $T$ and $\alpha\in (0,1)$, then $T\in \mathcal{A}(\alpha) $ is equivalent with any following statements:

(i): $||Tx-Ty ||^2 \leq ||x-y ||^2 -\frac{1-\alpha}{\alpha} || (x-y)-(Tx-Ty) ||^2, \forall x, y \in \Omega$.

(ii): $||Tx-Ty ||^2+ (1-2\alpha)||x-y||^2 \leq 2(1-\alpha) \langle x-y, Tx-Ty \rangle, \forall x, y \in \Omega$.
\end{lem}

\begin{lem}[Robbins-Siegmund]\label{lem_sequence}
Suppose nonnegative sequences $\{\alpha_k\}$, $\{\beta_k\}$ and $\{v_k\}$ satisfy the  recursive relations
$\alpha_{k+1}\leq \alpha_k -\beta_k+v_k,\forall k$
 and $\sum_{k=1}^{\infty}v_k <\infty$, then  $\{\alpha_k\}$ converges, $\sum_{k=1}^{\infty}\beta_k <\infty$ and $\lim_{k\rightarrow \infty}\beta_k=0$.
\end{lem}

\vskip 3mm
{\bf Proof of Theorem \ref{thm_pppa_convergence}}:
$\Phi(\varpi_{k}-\hat{\varpi}_k) \in  {\mathfrak{M}} \hat{\varpi}_k$ implies that $ \exists u_k \in {\mathfrak{M}} \hat{\varpi}_k $ such that $\Phi(\varpi_{k}-\hat{\varpi}_k) =u_k$.

Since $\Phi$ is positive definite, $\varpi_{k}-\hat{\varpi}_k=\Phi^{-1}u_k$. That is $\varpi_{k}-\hat{\varpi}_k\in \Phi^{-1}{\mathfrak{M}}\hat{\varpi}_k$. Since $\Phi$ is positive definite and ${\mathfrak{M}} $ is maximally monotone, we have $\Phi^{-1}{\mathfrak{M}}$ is  maximally monotone under the $\Phi-$induced norm $||\cdot||_{\Phi}$.
In fact, $\Phi$ is positive definite and nonsingular.
For any $(x,u)\in gra \Phi^{-1}{\mathfrak{M}}$ and $(y,v)\in gra \Phi^{-1}{\mathfrak{M}}$,
$\Phi u \in \Phi \Phi^{-1}{\mathfrak{M}}(x)\in {\mathfrak{M}}(x) $ and
$\Phi v \in \Phi \Phi^{-1}{\mathfrak{M}}(y)\in {\mathfrak{M}}(y) $.
Then $\langle x-y,u-v\rangle_{\Phi}= \langle  x-y, \Phi(u-v) \rangle \geq 0, \forall x,y\in dom \mathfrak{M}$,  since ${\mathfrak{M}}$ is monotone.   Therefore, $\Phi^{-1}{\mathfrak{M}}$ is monotone under the $\Phi-$induced inner product $\langle \cdot, \cdot \rangle_{\Phi}$. 
Furthermore, take $(y,v)$ with $y\in dom \mathfrak{M}$, and
$ \langle x-y, u-v \rangle_{\Phi} \geq 0$, for any other $(x,u) \in  gra(\Phi^{-1}{\mathfrak{M}})$.
For any $(x, \tilde{u}) \in gra {\mathfrak{M}}$, we have $(x,\Phi^{-1}\tilde{u})\in gra(\Phi^{-1}{\mathfrak{M}})$.
$ \langle  x-y, \Phi( \Phi^{-1}\tilde{u}-v) \rangle  \geq 0$, or equivalently,
$ \langle  x-y, \tilde{u}-\Phi v) \rangle  \geq 0$.
Since ${\mathfrak{M}}$ is maximally monotone, then $(y,\Phi v) \in gra \bar{\mathfrak{M}}$. We conclude that $v\in \Phi^{-1}\bar{\mathfrak{B}} (y)$ which implies that $\Phi^{-1}{\mathfrak{M}}$ is maximally monotone under $|| \cdot||_{\Phi}$. In the later proof, we will use $|| \cdot||$ for $||\cdot ||_{\Phi}$.

Therefore,  $\hat{\varpi}_k=({\rm Id}+\Phi^{-1}{\mathfrak{M}} )^{-1}\varpi_k$.
Denote $T=({\rm Id}+\Phi^{-1}{\mathfrak{M}} )^{-1}$, then
$T$ a firmly nonexpansive operator by Lemma \ref{lem_maxiamlly_monotone}. In other words,
their exists a nonexpansive operator $T^{'}$ such that $T=\frac{1}{2}{\rm Id}+\frac{1}{2}T^{'}$.
Hence $\hat{\varpi}_k=T\varpi_k= \frac{1}{2}\varpi_k+ \frac{1}{2}T^{'} \varpi_k$.
Moreover, given any $\varpi^* \in zer\Phi^{-1}{\mathfrak{M}}$, or equivalently, $\varpi^* \in zer{\mathfrak{M}}$, $\varpi^*$ is a fixed point of $T$ and $T^{'}$, i.e.,
$T\varpi^*=\varpi^*$ and $T^{'}\varpi^*=\varpi^*$,  with the definition of resolvent.

Denote $\breve{\varpi}_{k+1}=\varpi_k+\rho(\hat{\varpi}_k- \varpi_k )$. We have $
\breve{\varpi}_{k+1} = \varpi_k + \rho(\frac{1}{2}\varpi_k+\frac{1}{2}T^{'}\varpi_k - \varpi_k)
=(1- \frac{\rho}{2})\varpi_k+\frac{\rho}{2}T^{'}\varpi_k
$.
Denote $\tilde{T}= (1-\frac{\rho}{2}){\rm Id}+\frac{\rho}{2}T^{'}$, then $\tilde{T} \in \mathcal{A}(\frac{\rho}{2})$ since $\rho \in [1,2)$.
Moreover, given any $\varpi^*\in zer{\mathfrak{M}}$ we have $\tilde{T}\varpi^*= (1-\frac{\rho}{2})\varpi^*+\frac{\rho}{2}T^{'}\varpi^*=\varpi^*$ since $\varpi^*$ is a fixed point of $T^{'}$.

Given any $\varpi^*\in zer{\mathfrak{M}}$, with (i) of Lemma \eqref{lem_alpha_opertor} we have,
\begin{equation}
\begin{array}{l}\label{equ_thm_conver_1}
|| \breve{\varpi}_{k+1}- \varpi^*||^2
=|| \tilde{T}\varpi_k - \tilde{T}\varpi^*  ||^2 \\
\leq ||\varpi_k-\varpi^* ||^2 -\frac{2-\rho}{\rho} || \varpi_k-\varpi^*-(\tilde{T}\varpi_k-\tilde{T}\varpi^*) ||^2\\
= ||\varpi_k-\varpi^* ||^2
\end{array}
\end{equation}
Therefore, $|| \breve{\varpi}_{k+1}- \varpi^*|| \leq ||\varpi_k-\varpi^* ||$.
We also have $||\breve{\varpi}_{k+1}-\varpi_{k+1} ||\leq \rho \nu_k$ since $\breve{\varpi}_{k+1}-\varpi_{k+1}=\rho(\hat{\varpi}_k-\tilde{\varpi}_k)$ and $||\hat{\varpi}_k-\tilde{\varpi}_k ||\leq \nu_k$ due to Algorithm \ref{alg_pppa}. Then by the triangle inequality
\begin{equation}
\begin{array}{l} \nonumber
||{\varpi}_{k+1}-\varpi^* ||\leq || \breve{\varpi}_{k+1}-\varpi_{k+1} ||+|| \breve{\varpi}_{k+1}-\varpi^* ||\\
\leq ||\varpi_k-\varpi^* ||+ \rho \nu_k
\end{array}
\end{equation}

Since $\sum_{k=1}^{\infty} \rho \nu_k<\infty$, we conclude that $\{||{\varpi}_{k}-\varpi^* ||\}$ converges for any given $\varpi^* \in zer{\mathfrak{M}}$ with Lemma \ref{lem_sequence}.  Hence, $\{|| \breve{\varpi}_{k}- \varpi^*||\}$ and $\{||\varpi_k-\varpi^* ||\}$ are both bounded sequences, and we denote $c_4=\sup_{k} || \breve{\varpi}_{k}- \varpi^*||$.

Since $T=\frac{1}{2}{\rm Id}+\frac{1}{2}T^{'}$ is firmly nonexpansive, ${\rm Id}-T = \frac{1}{2} {\rm Id}+ \frac{1}{2}(-T^{'})$ is also firmly nonexpansive.
By (ii) of Lemma \ref{lem_alpha_opertor}, ${\rm Id}-T \in \mathcal{A}(\frac{1}{2})$ if and only if $\forall \varpi_1,\varpi_2 \in domT,$
\begin{equation}\label{firmly_nonexpansively}
\begin{array}{l}
||({\rm Id}-T)\varpi_1-({\rm Id}-T)\varpi_2 ||^2 \\
\leq \langle\varpi_1-\varpi_2, ({\rm Id}-T)\varpi_1-({\rm Id}-T)\varpi_2 \rangle
\end{array}
\end{equation}
Hence, we have
\begin{equation}
\begin{array}{l}\nonumber
|| \breve{\varpi}_{k+1}  -  \varpi^*||^2 =||\varpi_k+\rho(T\varpi_k- \varpi_k )-\varpi^* ||^2\\
=||\varpi_k-\rho({\rm Id}-T)\varpi_k-\varpi^* ||^2\\
=|| \varpi_k-\varpi^* ||^2 +  \rho^2||({\rm Id}-T)\varpi_k ||^2 \\
-  2\rho\langle \varpi_k-\varpi^*,({\rm Id}-T)\varpi_k- ({\rm Id}-T)\varpi^*\rangle\\
\leq || \varpi_k-\varpi^* ||^2 -(2\rho-\rho^2)||({\rm Id}-T)\varpi_k ||^2 \\
\end{array}
\end{equation}
where the third equality follows from $({\rm Id}-T)\varpi^*=\mathbf{0}$ and the last inequality follows from  \eqref{firmly_nonexpansively}. Denote $c_6=(2\rho-\rho^2)$, then we also have  $||\breve{\varpi}_{k+1}-\varpi_{k+1} ||^2\leq \rho^2 \nu^2_k$ and
\begin{equation}
\begin{array}{l}\nonumber
|| {\varpi}_{k+1}- \varpi^*||^2 = ||\varpi_{k+1}-\breve{\varpi}_{k+1}+\breve{\varpi}_{k+1}-\varpi^* ||^2\\
\leq \rho^2\nu_k^2 + ||\breve{\varpi}_{k+1}-\varpi^* ||^2+
2c_4\rho \nu_k\\
\leq ||\varpi_k-\varpi^* ||^2 -c_6
|| \varpi_k-T\varpi_k ||^2+ \rho (\rho\nu_k+2c_4)\nu_k
\end{array}
\end{equation}

We have  $\sum_{k=1}^{\infty}(\rho^2\nu_k+2\rho c_4)\nu_k <\infty$ due to $\sum_{k=1}^{\infty} \nu_k <\infty$.
By Lemma \ref{lem_sequence}, we conclude that  $\sum_{k=1}^{\infty}
|| \varpi_k-T\varpi_k ||^2 < \infty$, and $\lim_{k\rightarrow \infty} \varpi_k-T\varpi_k =\mathbf{0}$.

Since $\{|| \varpi_{k}-  \varpi^* || \}$ converges, $\{\varpi_k\}$ is a bounded sequence. There exists a subsequence$\{\varpi_{n_k}\}$  that converges to $\acute{\varpi}^*$.
Passing to limiting point of  Algorithm \ref{alg_pppa}, we have $T\acute{\varpi}^*=\acute{\varpi}^*$ by $\lim_{n_{k}\rightarrow \infty}T {\varpi}_{n_k}-\varpi_{n_k}=\mathbf{0}$ and (Lipschitz) continuity of $T$. Therefore, the limiting point $\acute{\varpi}^*$  is a fixed point of $T$ and is a zero of $\mathfrak{M}$ in \eqref{MI_agumented}.
Setting $\varpi^*=\acute{\varpi}^*$ in \eqref{equ_thm_conver_1}, we have  $\{|| \varpi_{k}-  \acute{\varpi}^* || \}$ is bounded   and converges. Since there exists a subsequence $\{\varpi_{n_k}\}$ that converges to $\acute{\varpi}^*$, it follows that  $\{|| \varpi_{k}-  \acute{\varpi}^* || \}$ converges to zero. Therefore, the whole sequence $\{\varpi_k\}$ generated from  Algorithm \ref{alg_pppa} with any initial point converges to  $\varpi^*$, and $\varpi^* \in zer{\mathfrak{M}}$.
\hfill $\Box$

\vskip 3mm

{\bf Proof of Theorem \ref{thm_equvalent}:}

We first give some useful relations derived from Algorithm \ref{alg_pppa} when $\mathfrak{M}=\bar{\mathfrak{M}}^e$ and $\Phi=\Phi^e$.

Write  $\Phi^e(\varpi_{k}-\hat{\varpi}_k) \in  \bar{\mathfrak{M}}^e \hat{\varpi}_k$ in its componentwise form
\begin{equation}
\begin{array}{l}\label{equ_thm_eqvalent_5}
R(\mathbf{x}^{'}_k-\hat{\mathbf{x}}^{'}_{k}) - \Lambda^{T}(\eta_k-\hat{\eta}_{k})+\Lambda^{T}(\theta_k-\hat{\theta}_{k}) \\
\qquad \in (N_{\bar{\Omega}}+ \partial F)(\hat{\mathbf{x}}^{'}_{k}) + \Lambda^T \hat{\eta}_{k}-\Lambda^T \hat{\theta}_{k}.\\
-\Lambda(\mathbf{x}^{'}_k-\hat{\mathbf{x}}^{'}_{k}) + 2H^{-1}(\eta_k-\hat{\eta}_{k})+\bar{V}(\mathbf{Z}^{'}_k-\hat{\mathbf{Z}}^{'}_{k}) \\
\qquad= -\Lambda \hat{\mathbf{x}}^{'}_{k}-\bar{V}\hat{\mathbf{Z}}^{'}_{k}+\bar{b}.\\
\bar{V}^T(\eta_k-\hat{\eta}_{k})+ W^{-1}(\mathbf{Z}^{'}_k-\hat{\mathbf{Z}}^{'}_{k})+\bar{V}^T(\theta_k-\hat{\theta}_{k}) \\
\qquad= \bar{V}^T \hat{\eta}_{k}-\bar{V}^T \hat{\theta}_{k}.\\
\Lambda(\mathbf{x}^{'}_k-\hat{\mathbf{x}}^{'}_{k}) + 2H^{-1}(\theta_k-\hat{\theta}_{k})+\bar{V}(\mathbf{Z}^{'}_k-\hat{\mathbf{Z}}^{'}_{k}) \\
\qquad= \Lambda \hat{\mathbf{x}}^{'}_{k}+\bar{V}\hat{\mathbf{Z}}^{'}_{k}-\bar{b}.
\end{array}
\end{equation}
Since $R,H$ and $W$ are positive definite,  \eqref{equ_thm_eqvalent_5}  gives 
\begin{align}
&R(\mathbf{x}^{'}_k-\hat{\mathbf{x}}^{'}_{k}) \in (N_{\bar{\Omega}}+ \partial F)(\hat{\mathbf{x}}^{'}_{k})
+\Lambda^T (\eta_{k}- \theta_{k})\label{equ_thm_eqvalent_1}.\\
&\hat{\eta}_{k}=\eta_k+\frac{H}{2}[\bar{V}\mathbf{Z}^{'}_k-\Lambda(\mathbf{x}^{'}_k-2\hat{\mathbf{x}}^{'}_{k})-\bar{b}]
\label{equ_thm_eqvalent_2}.\\
&\hat{\mathbf{Z}}^{'}_{k}=\mathbf{Z}^{'}_k+W\bar{V}^T(\eta_k-2\hat{\eta}_{k}+\theta_k).\label{equ_thm_eqvalent_3}\\
&\hat{\theta}_{k}=\theta_k+ \frac{H}{2}[ \Lambda(\mathbf{x}^{'}_k-2\hat{\mathbf{x}}^{'}_{k})  +\bar{V}(\mathbf{Z}^{'}_k-2\hat{\mathbf{Z}}^{'}_{k}) +\bar{b}].\label{equ_thm_eqvalent_4}
\end{align}
Denote
\begin{equation}\label{equ_thm_eqvalent_6}
\bar{\lambda}^{'}_k = \eta_k-\theta_k-H(\Lambda \mathbf{x}^{'}_k+\bar{V}\mathbf{Z}^{'}_k-\bar{b}),
\end{equation}
then by \eqref{equ_thm_eqvalent_1}
\begin{equation}
\begin{array}{l}\label{equ_thm_eqv_x}
R\mathbf{x}^{'}_k-\Lambda^T (\bar{\lambda}^{'}_k+H(\Lambda \mathbf{x}^{'}_k+\bar{V}\mathbf{Z}^{'}_k-\bar{b})) \in (N_{\bar{\Omega}}+ \partial F+R)\hat{\mathbf{x}}^{'}_{k}.
\end{array}
\end{equation}
Denote
$
\hat{\bar{\lambda}}^{'}_{k}=\hat{\eta}_{k}-\hat{\theta}_{k}-H(\Lambda \hat{\mathbf{x}}^{'}_{k}+\bar{V}\hat{\mathbf{Z}}^{'}_{k}-\bar{b}).
$
By \eqref{equ_thm_eqvalent_2} and \eqref{equ_thm_eqvalent_4}, we have
\begin{equation}
\begin{array}{l}\label{equ_thm_eqv_lambda}
\hat{\bar{\lambda}}^{'}_{k}=\eta_k+\frac{H}{2}[\bar{V}\mathbf{Z}^{'}_k-\Lambda(\mathbf{x}^{'}_k-2\hat{\mathbf{x}}^{'}_{k})-\bar{b}]\\
\quad -(\theta_k+ \frac{H}{2}[ \Lambda(\mathbf{x}^{'}_k-2\hat{\mathbf{x}}^{'}_{k})+\bar{V}(\mathbf{Z}^{'}_k-2\hat{\mathbf{Z}}^{'}_{k}) +\bar{b}])\\
\quad -H(\Lambda \hat{\mathbf{x}}^{'}_{k}+\bar{V}\hat{\mathbf{Z}}^{'}_{k}-\bar{b})\\
=\eta_k-\theta_k+H[-\Lambda(\mathbf{x}^{'}_k-\hat{\mathbf{x}}^{'}_{k})]\\
=\bar{\lambda}^{'}_k+ H(\Lambda \mathbf{x}^{'}_k+\bar{V}\mathbf{Z}^{'}_k-\bar{b})+H[-\Lambda(\mathbf{x}^{'}_k-\hat{\mathbf{x}}^{'}_{k})]\\
=\bar{\lambda}^{'}_k+H(\Lambda \hat{\mathbf{x}}^{'}_{k}+\bar{V}\mathbf{Z}^{'}_k-\bar{b}).
\end{array}
\end{equation}
From \eqref{equ_thm_eqvalent_6}, \eqref{equ_thm_eqv_lambda}, $\hat{\bar{\lambda}}^{'}_{k}=\eta_k-\theta_k+
H[-\Lambda(\mathbf{x}^{'}_k-\hat{\mathbf{x}}^{'}_{k})]$.

Then by \eqref{equ_thm_eqvalent_3} and \eqref{equ_thm_eqvalent_2}
\begin{equation}
\begin{array}{lll}\label{equ_thm_eqv_z}
\hat{\mathbf{Z}}^{'}_{k}=\mathbf{Z}^{'}_k+W\bar{V}^T(\theta_k-\eta_k
-{H}[\bar{V}\mathbf{Z}^{'}_k-\Lambda(\mathbf{x}^{'}_k-2\hat{\mathbf{x}}^{'}_{k})-\bar{b}] ) \\
=\mathbf{Z}^{'}_k+W\bar{V}^T(-\hat{\bar{\lambda}}^{'}_{k}+H[-\Lambda(\mathbf{x}^{'}_k-\hat{\mathbf{x}}^{'}_{k})]\\
\qquad-{H}[\bar{V}\mathbf{Z}^{'}_k-\Lambda(\mathbf{x}^{'}_k-2\hat{\mathbf{x}}^{'}_{k})-\bar{b}] ) \\
=\mathbf{Z}^{'}_k-W\bar{V}^T(\hat{\bar{\lambda}}^{'}_{k}+H[\Lambda\hat{\mathbf{x}}^{'}_{k}+
\bar{V}\mathbf{Z}^{'}_k-\bar{b}] ).
\end{array}
\end{equation}

\vskip 1mm
Then we prove \eqref{correspondence} by induction.
Firstly, we choose
$\mathbf{x}^{'}_{0}=\mathbf{x}_{0}, \eta_0 = \bar{\lambda}_0+H(\Lambda \mathbf{x}_{0}+\bar{V}^T \bar{\lambda}_0-\bar{b} ),\mathbf{Z}^{'}_0=\mathbf{Z}_0,\theta_0=\mathbf{0}$.
Hence \eqref{correspondence} holds at $k=0$.


Suppose \eqref{correspondence} is true at time $k$, then $\bar{\lambda}^{'}_k$ defined in \eqref{equ_thm_eqvalent_6}
has $\bar{\lambda}^{'}_k=\bar{\lambda}_k$.
We can choose $\hat{\mathbf{x}}^{'}_k = \hat{\mathbf{x}}_k$ to satisfy equation \eqref{equ_thm_eqv_x} due to $\mathbf{x}^{'}_k=\mathbf{x}_k$, $\mathbf{Z}_k^{'}=\mathbf{Z}_k$ and \eqref{alg_compact}.
Then we choose $\tilde{\mathbf{x}}^{'}_k = \tilde{\mathbf{x}}_k$ such that
$||\hat{\mathbf{x}}^{'}_k- \tilde{\mathbf{x}}^{'}_k||=||\hat{\mathbf{x}}_k-\tilde{\mathbf{x}}_k || \leq \mu_k$.
Thereby, we have $\mathbf{x}^{'}_{k+1}= \mathbf{x}^{'}_k+\rho( \tilde{\mathbf{x}}^{'}_k- \mathbf{x}^{'}_k)=\mathbf{x}_{k+1}$.

Recall that $\hat{\eta}_k$ and $\hat{\theta}_k$ are generated by \eqref{equ_thm_eqvalent_2} and \eqref{equ_thm_eqvalent_4} from $\eta_k$ and $\theta_k$. Due to \eqref{equ_thm_eqv_lambda}, \eqref{correspondence} and   $\hat{\mathbf{x}}^{'}_k = \hat{\mathbf{x}}_k$ we have: 
\begin{equation}
\begin{array}{l}\label{equ_thm_eqvalent_7}
\hat{\bar{\lambda}}^{'}_{k}=\bar{\lambda}^{'}_k+H(\Lambda \hat{\mathbf{x}}^{'}_{k}+\bar{V}\mathbf{Z}^{'}_k-\bar{b})\\
\quad   \;              =\bar{\lambda}_k+H(\Lambda \hat{\mathbf{x}}_{k}+\bar{V}\mathbf{Z}_k-\bar{b}).
\end{array}
\end{equation}
By  \eqref{alg_compact} we also have $\tilde{\bar{\lambda}}_k =\bar{\lambda}_k+H(\Lambda \tilde{\mathbf{x}}_{k}+\bar{V}\mathbf{Z}_k-\bar{b}) $. Hence $\tilde{\bar{\lambda}}_k -\hat{\bar{\lambda}}^{'}_{k}= H\Lambda(\tilde{\mathbf{x}}_{k}-\hat{\mathbf{x}}_{k})$.
By \eqref{equ_thm_eqv_z}
 \begin{equation}
\begin{array}{l}
\hat{\mathbf{Z}}^{'}_{k}
=\mathbf{Z}^{'}_k-W\bar{V}^T(\hat{\bar{\lambda}}^{'}_{k}+H[\Lambda\hat{\mathbf{x}}^{'}_{k}+
\bar{V}\mathbf{Z}^{'}_k-\bar{b}] )\\
\qquad =\mathbf{Z}_k-W\bar{V}^T(\hat{\bar{\lambda}}^{'}_{k}+H[\Lambda\hat{\mathbf{x}}_{k}+
\bar{V}\mathbf{Z}_k-\bar{b}] ).
\end{array}
\end{equation}
We  choose $ \tilde{\mathbf{Z}}^{'}_{k}=\tilde{\mathbf{Z}}_{k}$ where
$\tilde{\mathbf{Z}}_{k}=\mathbf{Z}_k-W\bar{V}^T(\tilde{\bar{\lambda}}_{k}
+H(\Lambda \tilde{\mathbf{x}}_{k}+\bar{V}\mathbf{Z}_k-\bar{b})$ due to \eqref{alg_compact}, so that
\begin{equation}
\begin{array}{l}
\tilde{\mathbf{Z}}^{'}_{k}-\hat{\mathbf{Z}}^{'}_{k}= \tilde{\mathbf{Z}}_{k}-\hat{\mathbf{Z}}^{'}_{k}\\
=-W\bar{V}^T(\tilde{\bar{\lambda}}_{k}-\hat{\bar{\lambda}}^{'}_{k})-W\bar{V}^TH\Lambda(\tilde{\mathbf{x}}_{k}-\hat{\mathbf{x}}_{k})\\
=-2W\bar{V}^TH\Lambda(\tilde{\mathbf{x}}_{k}-\hat{\mathbf{x}}_{k}).
\end{array}
\end{equation}
Therefore, $\exists c_1>0$, such that $|| \tilde{\mathbf{Z}}^{'}_{k}-\hat{\mathbf{Z}}^{'}_{k}||\leq c_1\mu_k$.
Since  $\tilde{\mathbf{Z}}^{'}_{k}=\tilde{\mathbf{Z}}_{k}$ and $\mathbf{Z}^{'}_k=\mathbf{Z}_k$,
we have $\mathbf{Z}^{'}_{k+1} = \mathbf{Z}^{'}_k+\rho(\tilde{\mathbf{Z}}^{'}_k- \mathbf{Z}^{'}_k)=\mathbf{Z}^{k+1}$.

 Denote $\tilde{\bar{\lambda}}^{'}_k = \tilde{\eta}_{k}-\tilde{\theta}_{k}-H(\Lambda \tilde{\mathbf{x}}^{'}_{k}+\bar{V}\tilde{\mathbf{Z}}^{'}_{k}-\bar{b})
$. Then we want to find $\tilde{\eta}_k, \tilde{\theta}_k$ with $||\tilde{\eta}_k-\hat{\eta}_k||\leq c_2\mu_k$, $||\tilde{\theta}_k-\hat{\theta}_k||\leq c_3 \mu_k$
such that $\tilde{\bar{\lambda}}_k=\tilde{\bar{\lambda}}^{'}_k$.
Suppose $\tilde{\theta}_k$ and $\tilde{\eta}_k$ are chosen to ensure
$\tilde{\eta}_{k}-\tilde{\theta}_{k}-H(\Lambda \tilde{\mathbf{x}}^{'}_{k}+\bar{V}\tilde{\mathbf{Z}}^{'}_{k}-\bar{b})
 =\hat{\eta}_{k}-\hat{\theta}_{k}-H(\Lambda \hat{\mathbf{x}}^{'}_{k}+\bar{V}\hat{\mathbf{Z}}^{'}_{k}-\bar{b})+H\Lambda(\tilde{\mathbf{x}}_{k}-\hat{\mathbf{x}}_{k})$. Then  due to \eqref{alg_compact} and \eqref{equ_thm_eqvalent_7},
\begin{equation}
\begin{array}{l}\label{equ_thm_eqvalent_8}
\tilde{\bar{\lambda}}_k  =\bar{\lambda}_k + H(\Lambda \tilde{\mathbf{x}}_{k}+\bar{V}\mathbf{Z}_k-\bar{b})\\
=\hat{\bar{\lambda}}^{'}_{k}-H(\Lambda \hat{\mathbf{x}}_{k}+\bar{V}\mathbf{Z}_k-\bar{b})+ H(\Lambda \tilde{\mathbf{x}}_{k}+\bar{V}\mathbf{Z}_k-\bar{b})\\
=\hat{\bar{\lambda}}^{'}_{k}+H\Lambda(\tilde{\mathbf{x}}_{k}-\hat{\mathbf{x}}_{k})\\
=\hat{\eta}_{k}-\hat{\theta}_{k}-H(\Lambda \hat{\mathbf{x}}^{'}_{k}+\bar{V}\hat{\mathbf{Z}}^{'}_{k}-\bar{b})+H\Lambda(\tilde{\mathbf{x}}_{k}-\hat{\mathbf{x}}_{k})\\
=\tilde{\eta}_{k}-\tilde{\theta}_{k}-H(\Lambda \tilde{\mathbf{x}}^{'}_{k}+\bar{V}\tilde{\mathbf{Z}}^{'}_{k}-\bar{b})=\tilde{\bar{\lambda}}^{'}_k. \\
\end{array}
\end{equation}
Hence, let $\tilde{\theta}_k$ and $\tilde{\eta}_k$ be chosen as,
\begin{equation}
\begin{array}{l}
\tilde{\eta}_k=\hat{\eta}_k+H\Lambda(\tilde{\mathbf{x}}_{k}-\hat{\mathbf{x}}_{k})+\frac{1}{2} H\bar{V}(\tilde{\mathbf{Z}}^{'}_{k}-\hat{\mathbf{Z}}^{'}_{k}). \\
\tilde{\theta}_k=\hat{\theta}_k- H\Lambda(\tilde{\mathbf{x}}_{k}-\hat{\mathbf{x}}_{k})-\frac{1}{2} H\bar{V}(\tilde{\mathbf{Z}}^{'}_{k}-\hat{\mathbf{Z}}^{'}_{k}).
\end{array}
\end{equation}
Obviously, in this case $\exists c_2>0, c_3>0$ such that $||\tilde{\eta}_k-\hat{\eta}_k||\leq c_2\mu_k$, $||\tilde{\theta}_k-\hat{\theta}_k||\leq c_3 \mu_k$. Moreover, from \eqref{equ_thm_eqvalent_8} we have $\tilde{\bar{\lambda}}^{'}_k=\tilde{\bar{\lambda}}_k$.
Hence, we obtain that
\begin{equation}
\begin{array}{ll}\nonumber
\bar{\lambda}^{'}_{k+1}  &=  \eta_{k+1}-\theta_{k+1}-H(\Lambda \mathbf{x}^{'}_{k+1}+\bar{V}\mathbf{Z}^{'}_{k+1}-\bar{b})\\
\quad     & =(1-\rho)\eta_{k}+\rho \tilde{\eta}_k-[(1-\rho)\theta_{k}+\rho \tilde{\theta}_k]\\
\quad     &-H(\Lambda [ (1-\rho)\mathbf{x}^{'}_{k}+\rho \tilde{\mathbf{x}}^{'}_k]+\bar{V}[(1-\rho)\mathbf{Z}^{'}_{k}+\rho \tilde{\mathbf{Z}}^{'}_k]-\bar{b})\\
\qquad &=(1-\rho)[\eta_{k}-\theta_{k}-H(\Lambda\mathbf{x}^{'}_{k}+\bar{V}\mathbf{Z}^{'}_{k}-\bar{b})]\\
\quad &+\rho[ \tilde{\eta}_k -\tilde{\theta}_k-H(\Lambda \tilde{\mathbf{x}}^{'}_k+\bar{V} \tilde{\mathbf{Z}}^{'}_k-\bar{b}) ]\\
\qquad &= (1-\rho)\bar{\lambda}_k+\rho \tilde{\bar{\lambda}}_k=\bar{\lambda}_{k+1}.
\end{array}
\end{equation}

Therefore, when \eqref{correspondence} holds at time $k$, it also holds at time $k+1$.
Thus,  we have shown by induction that
given  sequences $col(\mathbf{x}_k,\mathbf{Z}_k,\bar{\lambda}_k)$
generated from Algorithm \ref{alg_d}
with initial points $\mathbf{x}_0,\mathbf{Z}_0,\bar{\lambda}_0$ and $\{\mu_k\}$, we can find  sequences
$\varpi_k=col(\mathbf{x}^{'}_k,\eta_k,\mathbf{Z}^{'}_k,\theta_k)$
generated from Algorithm \ref{alg_pppa} with  $\nu_k \leq \sqrt{1+c_1^2+c_2^2+c_3^2}\mu_k$ such that \eqref{correspondence} holds. Since $\sum_{k=1}^{\infty} \mu_k < \infty$, we have $\sum_{k=1}^\infty \nu_k <\infty$, and the conclusion follows.
\hfill $\Box$

\end{document}